\documentclass[12pt,reqno]{amsart}
\numberwithin{equation}{section}

\allowdisplaybreaks

\usepackage[T1]{fontenc}
\usepackage[sc]{mathpazo}
\usepackage[centering,scale=0.80]{geometry}
\usepackage{eucal}
\usepackage[all]{xy}
\usepackage{tikz}
\usepackage{wrapfig}
\usepackage{graphicx}

\usepackage{xcolor}

\usepackage[unicode,colorlinks=true,linkcolor=blue!45!black,citecolor=blue!45!black,urlcolor=blue!45!black]{hyperref}

\allowdisplaybreaks[1]

\usepackage{amsmath,amssymb,amsthm,mathtools}

\def\HHom{\operatorname{\mathcal{H}om}}
\def\EEnd{\operatorname{\mathcal{E}nd}}
\def\Hom{\operatorname{Hom}}
\def\End{\operatorname{End}}
\def\Ext{\operatorname{Ext}}
\def\Spec{\operatorname{Spec}}

\def\OO{\mathcal{O}}

\newtheorem*{Thm}{Theorem}
\newtheorem{maintheorem}{Theorem}

\newtheorem{theorem}{Theorem}[section]
\newtheorem{proposition}[theorem]{Proposition}
\newtheorem{lemma}[theorem]{Lemma}
\newtheorem{corollary}[theorem]{Corollary}
\theoremstyle{definition}

\theoremstyle{remark}
\newtheorem{remark}[theorem]{Remark}

\begin{document}

\title{Macdonald intersections and the $N!/k$ conjecture}

\author{Peter L. Guo}
\address[Peter L. Guo]{Center for Combinatorics, LPMC, 
Nankai University, Tianjin 300071, P.R. China}
\email{lguo@nankai.edu.cn}

\author{Mingyang Kang}
\address[Mingyang Kang]{Center for Combinatorics, LPMC, 
Nankai University, Tianjin 300071, P.R. China}
\email{2120240005@mail.nankai.edu.cn}

\author{Changjian Su}
\address[Changjian Su]{Yau Mathematical Sciences Center, Tsinghua University, Beijing 100084, P.R. China}
\email{changjiansu@mail.tsinghua.edu.cn}

\author{Rui Xiong}
\address[Rui Xiong]{Department of Mathematics, University of Michigan, Ann Arbor, MI 48109}
\email{xiongrui@umich.edu}

\begin{abstract}
We prove the Frobenius characteristic formula for the intersection of Garsia--Haiman modules predicted in the Science Fiction of  Bergeron and Garsia. 
As a consequence, we resolve the $N!/k$ conjecture. The proof uses  the geometry of the  Hilbert scheme of points in the plane developed by Haiman. 
\end{abstract}



\maketitle

\tableofcontents

\section{Introduction}

Let $q,t$ be independent parameters, and let $X=x_1+x_2+\cdots$ be an
alphabet of independent variables. Write $s_\lambda[X]$ for the Schur function
indexed by a partition $\lambda$.
The \emph{modified Macdonald polynomials}
$\widetilde H_\lambda[X;q,t]$ are the transformed integral forms of the Macdonald polynomials \cite{Mac-2,Mac-3}.
These polynomials are homogeneous symmetric functions  over $\mathbb{Q}(q,t)$ characterized by 
\cite[Proposition 2.6]{Haiman-JJ} 
\begin{enumerate}
    \item 
    $\widetilde{H}_\lambda[(1-q)X;q,t]
\in\operatorname{span}\{s_\mu[X]:\mu\geq\lambda\}$; 
    \item 
$\widetilde{H}_\lambda[(1-t)X;q,t]
\in\operatorname{span}
\{s_\mu[X]:\mu\geq\lambda'\}$;
    \item
$\widetilde H_\lambda[1;q,t]=1$. 
\end{enumerate}
Here $\lambda'$ is the conjugate partition of $\lambda$, and $\geq$ denotes  the dominance partial order. 

The Macdonald positivity conjecture \cite{Mac-2}  asserts that 
\[
\widetilde{H}_\lambda[X;q,t]\in
\bigoplus_{\mu}\mathbb{Z}_{\geq 0}[q,t]s_\mu(X).
\] 
To attack this conjecture, 
 Garsia and Haiman   \cite{GH, GH-22} introduced 
for each  $\lambda\vdash n$ a bigraded module $\mathsf{GH}_\lambda$ of 
the symmetric group $S_n$:
$$
\mathsf{GH}_\lambda=
\mathbb C\bigl[
\tfrac{\partial}{\partial x_1},\tfrac{\partial}{\partial y_1},\ldots,
\tfrac{\partial}{\partial x_n},\tfrac{\partial}{\partial y_n}
\bigr]\Delta_\lambda,
\qquad
\Delta_\lambda=
\det(x_i^a y_i^b)_{\substack{1\leq i\leq n\\(a,b)\in\lambda}},
$$
where a partition is identified with   its diagram 
\[\{(j-1,i-1):1\leq i\leq \ell(\lambda), 1\leq j\leq \lambda_{i}\},\]
and conjectured  that    the bigraded Frobenius characteristic  $\operatorname{Frob}_{q,t}$ of $\mathsf{GH}_\lambda$ is $\widetilde{H}_\lambda[X;q,t]$.

In a major  breakthrough in algebraic combinatorics, 
Haiman \cite{Haiman} resolved   both the  conjecture stated above and the  $n!$ conjecture of Garsia and Haiman  \cite{GH}  by developing  fundamental geometric properties of the Hilbert scheme of points  in the
plane. 

\begin{Thm}[{Haiman \cite[Theorem 3.2]{Haiman}}]
For every partition $\lambda\vdash n$,
$$
\operatorname{Frob}_{q,t}\mathsf{GH}_\lambda
=\widetilde H_\lambda[X;q,t],
\qquad \dim_{\mathbb C}\mathsf{GH}_\lambda=n!.
$$ 
\end{Thm}

The ``Science Fiction''  in the seminal work  of  Bergeron
and Garsia  \cite{BG} proposed  a package of ``heuristics'' to investigate  the     Garsia--Haiman submodules indexed by 
partitions obtained by removing  corners from  a given partition $\lambda\vdash n$. The modern target of the Science Fiction centers around   predicted  formulas for the  Frobenius characteristic and the dimension of the intersection of these submodules.

The main goal of this paper is to prove the above  predictions  in the Science Fiction.

\begin{Thm}[{=Theorem \ref{thm:Section6}, Prediction in \cite[(3.26)]{BG}}]
Let $n\geq 2$ and $k\geq1$ be integers, and denote $N=n-1$. 
For a given partition  
$\lambda\vdash n$, let $\mu(1),\ldots,\mu(k)\vdash N$ be distinct
partitions each of which is obtained from $\lambda$ by removing a corner. Then
\begin{equation}\label{eq:main}
\operatorname{Frob}_{q,t}\left(
\mathsf{GH}_{\mu(1)}\cap \cdots \cap 
\mathsf{GH}_{\mu(k)}\right)
=\sum_{i=1}^k
\left(\prod_{\substack{1\leq j\leq k\\j\ne i}}
\frac{T_{\mu(j)}}{T_{\mu(j)}-T_{\mu(i)}}\right)
\widetilde H_{\mu(i)}[X;q,t],
\end{equation}
where $T_\rho=\prod_{(a,b)\in\rho}q^a t^b$ for a partition $\rho$.
\end{Thm}

As a significant application, according to   \cite[Equation (3.28)]{BG}, the establishment of   \eqref{eq:main} implies the following $N!/k$ conjecture. 

\begin{corollary}[$N!/k$ conjecture, {Prediction in \cite[(I.34)]{BG}}]\label{cor:N!/k}
Under the same hypotheses, we have
$$
\dim_{\mathbb C}\left(
\mathsf{GH}_{\mu(1)}\cap \cdots \cap 
\mathsf{GH}_{\mu(k)}\right)=\frac{N!}{k}.
$$
\end{corollary}

Before proceeding, we  review some prior progress. 
\begin{itemize}
    \item When the partition is a hook shape, Bergeron and Hamel \cite{BH} obtained a formula for
the intersection character, and Armon \cite{Armon} constructed an
explicit basis proving the dimension formula.

    \item   Kim, Lee and Oh \cite{KLO}
obtained shuffle and specialization formulas for the predicted characters. Kim and Oh \cite{KO}  studied a generalization of the   intersection character and  the $N!/k$ conjecture by defining  Macdonald
piece polynomials.

\item The  $k=2$ case has received considerable  attention. In this case,  Butler \cite{Butler1994} conjectured  that the right-hand side of \eqref{eq:main} is Schur positive.  
 Kim,   Lee and   Oh
\cite{KLO2026} proved that the character is monomial positive, and confirmed Butler's conjecture in some special situations. Butler's conjecture was recently settled in full by Guo, Kang and Xiong \cite{GKX}  
using the Procesi bundle on a projective line; see Khoroshkin and Makedonskyi
\cite{KM} for an alternative representation-theoretic proof. 

\item 
   Khoroshkin and Makedonskyi
\cite{KM-2}  proved the intersection character formula and the dimension formula for partitions with at most two rows by using generalized parabolic Demazure modules.

\end{itemize}

Our  proof builds on Haiman's geometric construction \cite{Haiman}, but requires    establishing new geometric properties that are essential to the argument, which will be briefly outlined.
The Hilbert scheme $H_n$ of $n$ points in $\mathbb{A}^2$ parametrizes ideals in $\mathbb{C}[x,y]$ of colength $n$. 
For each partition $\lambda\vdash n$, we have a monomial ideal $I_\lambda\in H_n$ generated by the monomials outside $\lambda$. 
The most important ingredient is the Procesi bundle $\mathcal{P}_n$, recalled in Section \ref{sec:prelimiaries}, putting polynomial quotient models of the Garsia--Haiman modules in one family. 
That is, the fiber of the Procesi bundle $\mathcal{P}_n$ at $I_\lambda$ is dual to the Garsia--Haiman module $\mathsf{GH}_\lambda$ in each bidegree.

In our case, we will not only consider the fiber at a single point. 
Instead, we connect the chosen corners
by a projective space. 
Let $V\subset \mathbb{C}[x,y]$ be the span of the monomials of the $k$ removable corners. 
We obtain a morphism 
$$
f:\mathbb P(V)\longrightarrow H_{n-1},
\qquad L\longmapsto I_\lambda+L.
$$
The first observation is that the right-hand side of \eqref{eq:main} can be identified with the Atiyah–Bott localization formula for the alternating sum
$$\sum_{i\geq 0}
(-1)^i\operatorname{Frob}_{q,t}
\big(H^i(\mathbb{P}(V),f^*\mathcal{P}_{n-1}^\vee)^\vee\big).$$
This motivates us to equip the intersection on the left-hand side of \eqref{eq:main} with a cohomological meaning, namely, 
it is dual to the space of global sections
$$
H^0\bigl(\mathbb P(V),f^*\mathcal P_{n-1}^{\vee}\bigr)^{\vee}
$$
in each bidegree. 
Moreover, we also need the vanishing of higher cohomology to ensure localization computes $H^0$. 
See Theorem \ref{thm:Section5} for the precise statement. 

The main difficulty is to justify this identification. 
In order to prove it, we need to know the positivity (in a proper sense) of the pullback of the Procesi bundle $f^*\mathcal{P}_{n-1}$. 
The Procesi bundle is globally generated, but that is not enough for our purpose. 
It turns out that  the Castelnuovo--Mumford regularity over projective space is the positivity we are looking for. 
We will establish this by showing that $f^*\mathcal{P}_{n-1}$ is a direct summand of direct sums of exterior powers of the universal quotient bundle, which is known to be $0$-regular; see Theorem \ref{thm:Section4}. 

The proof strategy is to lift the problem to the nested Hilbert scheme $H_{n-1,n}$, parameterizing pairs of nested ideals of colengths $n$ and $n-1$, respectively. 
The morphism $f$ above lifts to $\tilde{f}:\mathbb{P}(V)\to H_{n-1,n}$ by attaching $I_\lambda$ to the morphism $f$.  
A technique from tilting theory will reduce the problem to the vanishing of $\Ext$-groups. 
These will follow from
Theorem \ref{thm:Section3} stating that the relevant morphism is cohomologically trivial. 

A special case of Theorem~\ref{thm:Section3} was proved by Haiman \cite{Haiman} in order to prove
that the isospectral Hilbert scheme $X_n$ is Gorenstein. 
Our approach is carried out along similar lines.

The geometric argument also gives a simultaneous decomposition
of the Garsia--Haiman modules and their intersections, which responds    affirmatively   to the predictions \cite[(I.23) and (I.24)]{BG} in   the original  Science Fiction.

\begin{Thm}[=Theorem \ref{lem:primitivepro}]
Under the same hypotheses, write 
$\mathsf{GH}_A=\bigcap_{i\in A}\mathsf{GH}_{\mu(i)}$ for every  nonempty subset $A\subseteq [k]:=\{1,\ldots, k\}$. 
There exist bigraded $S_{n-1}$-submodules
$\mathsf{GH}'_A\subset \mathsf{GH}_A$, called 
\emph{primitive Garsia--Haiman modules}, such that  
\begin{equation}\label{eq:main2}
\mathsf{GH}_{\mu(1)}
+\cdots+\mathsf{GH}_{\mu(k)}
=\bigoplus_{\varnothing \neq A\subseteq [k] } \mathsf{GH}'_A,\qquad 
\mathsf{GH}_{A}=\bigoplus_{B\supseteq A} \mathsf{GH}'_B.
\end{equation}
\end{Thm}

\begin{figure}[h t]
$$\begin{tikzpicture}
\draw[fill = red, opacity = 0.15] ({1*cos(30)},{1*sin(30)}) circle (1.7);
\node at ({3.2*cos(30)},{3.2*sin(30)}) 
{\(\mathsf{GH}_{\mu(1)}\)};
\draw[fill = blue, opacity = 0.15] ({1*cos(150)},{1*sin(150)}) circle (1.7);
\node at ({3.3*cos(150)},{3.3*sin(150)}) 
{\(\mathsf{GH}_{\mu(2)}\)};
\draw[fill = yellow!70!black, opacity = 0.3] ({1*cos(270)},{1*sin(270)}) circle (1.7);
\node at ({3.0*cos(270)},{3.0*sin(270)}) 
{\(\mathsf{GH}_{\mu(3)}\)};
\node at ({1.8*cos(30)},{1.8*sin(30)}) 
{\(\mathsf{GH}'_{\{1\}}\)};
\node at ({1.8*cos(150)},{1.8*sin(150)}) 
{\(\mathsf{GH}'_{\{2\}}\)};
\node at ({1.8*cos(270)},{1.8*sin(270)}) 
{\(\mathsf{GH}'_{\{3\}}\)};
\node at ({1.1*cos(90)},{1.1*sin(90)}) 
{\(\scriptstyle\mathsf{GH}'_{\{1,2\}}\)};
\node at ({1.1*cos(210)},{1.1*sin(210)}) 
{\(\scriptstyle\mathsf{GH}'_{\{2,3\}}\)};
\node at ({1.1*cos(330)},{1.1*sin(330)}) 
{\(\scriptstyle\quad\mathsf{GH}'_{\{1,3\}}\)};
\node at (0,0) 
{\(\scriptstyle\mathsf{GH}'_{\{1,2,3\}}\)};
\end{tikzpicture}$$
    \caption{Primitive Garsia--Haiman modules}
    \label{fig:primGH}
\end{figure}
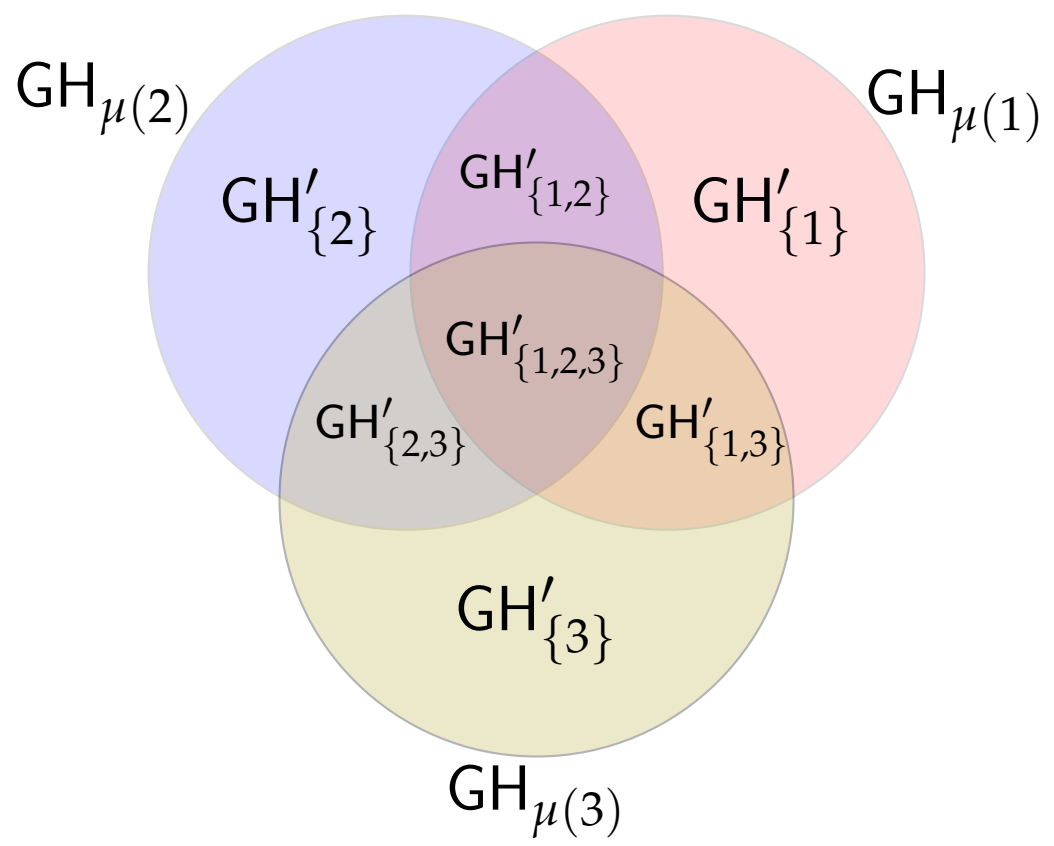

An illustration of the assertions in \eqref{eq:main2} is given in Figure \ref{fig:primGH}.

If we pick a homogeneous basis for each $\mathsf{GH}'_A$, then their union satisfies the properties required  in \cite[(I.23) and (I.24)]{BG}. 
In fact, both   primitive Garsia--Haiman modules and the sum of Garsia--Haiman modules in \eqref{eq:main2}  admit  cohomological meanings, see Section \ref{sec:localization}. 
Note also that the statement in \eqref{eq:main2} suggests   that $\mathsf{GH}_{\mu(1)},\ldots,\mathsf{GH}_{\mu(k)}$ possess  a distributive lattice structure in $\mathbb{C}[x_1,y_1,\ldots,x_{n-1},y_{n-1}]$.

The paper is organized as follows. Section~\ref{sec:prelimiaries} recalls Hilbert schemes,
Procesi bundles, and the local Hilbert--Burch equations. Section~\ref{sec:polygraph} proves
the mixed polygraph pushforward theorem. Section~\ref{sec:tilting} uses tilting theory to obtain
the local direct-summand statement. Section \ref{sec:regular} proves regularity and
identifies the common quotient with cohomology. Section \ref{sec:localization} applies
localization to prove \eqref{eq:main} and \eqref{eq:main2}.

\subsection*{Acknowledgments}
This work is supported by the National Key R\&D Program of China (No. 2024YFA1014700), the National Natural Science Foundation of China (No. 12371329), and the Fundamental Research Funds for the Central Universities (No. 63263094).

\subsection*{AI disclosure}

This article is a follow-up to our previous work~\cite{GKX} on Butler's conjecture. During its preparation, we used ChatGPT interactively to help check calculations and explore possible proof strategies. The observation linking the conjectural formulas \eqref{eq:main} and \eqref{eq:main2} to sheaf cohomology through equivariant localization is due entirely to the authors. All AI-generated material was reviewed and verified by the authors, who take full responsibility for the final content.

\section{Preliminaries: Hilbert schemes, Procesi bundles, etc.}\label{sec:prelimiaries}

In this section, we recall Haiman's theory of Hilbert schemes. 
We use the following diagrams:
$$
\begin{matrix}
\xymatrix@!=2pc{
X_n\ar[d]_{\rho_n}\ar[r] & (\mathbb{A}^2)^n\ar[d]\\
H_n\ar[r] & S^n(\mathbb{A}^2),}
\end{matrix}\qquad 
\begin{matrix}
\xymatrix@!=2pc{
U_n\ar[d]_{\pi_n}\ar[r] & \mathbb{A}^2\\
H_n& }
\end{matrix}\quad \text{ and } 
\begin{matrix}
\xymatrix@!=2pc{
& 
H_{n-1,n}\ar[dl]_{q_{n-1}}\ar[dr]_{p_n}
\ar[r]^-{\alpha}
& U_n\ar[d]^{\pi_n}
\\
H_{n-1} && H_n.}
\end{matrix}
$$

\subsection{Hilbert schemes}\label{2-11}

For $n\geq1$, let
$$H_n
=\{I\trianglelefteq \mathbb{C}[x,y]: 
\dim_{\mathbb{C}} \mathbb{C}[x,y]/I=n\}$$
be the \emph{Hilbert scheme} of $n$ points in  $\mathbb{A}^2$. 
By Fogarty's theorem \cite{Fogarty}, the scheme $H_n$ is a smooth irreducible variety of dimension $2n$. Let $S^n(\mathbb{A}^2)$ be the $n$-th symmetric power of $\mathbb{A}^2$. 
It parametrizes effective $0$-cycles on $\mathbb{A}^2$ of degree $n$. There is a Hilbert--Chow map from $H_n$ to $S^n(\mathbb{A}^2)$, recording the support (with multiplicities) of the finite subscheme $\Spec(\mathbb{C}[x,y]/I)$. 

Let $\pi_n:U_n\to H_n$ be the \emph{universal family}. Its closed points can be described as follows: 
$$U_n = \{(I,a): a\in \operatorname{Supp}(\Spec(\mathbb{C}[x,y]/I))\}\subset 
H_n\times \mathbb{A}^2.$$
It follows from the definition of the Hilbert scheme that $\pi_n$ is finite flat of degree $n$. 
The \emph{universal bundle} over $H_n$ is defined to be
$$\mathcal{B}_n = \pi_{n*}\OO_{U_n}.$$
That is, for $I\in H_n$, the fiber
$\mathcal{B}_n(I) = \mathbb{C}[x,y]/I$.

The \emph{isospectral Hilbert scheme} is the reduced fiber product
$X_n = \big(H_n\times_{S^n(\mathbb{A}^2)}(\mathbb{A}^2)^n\big)_{\rm red}$. 
That is, 
$$X_n = \{(I,P_1,\ldots,P_n): 
[P_1]+\ldots+[P_n] = [\mathbb{C}[x,y]/I]\}
\subset H_n\times (\mathbb{A}^2)^n,$$
where $[\mathbb{C}[x,y]/I]$ denotes the image of $I$ under the Hilbert--Chow morphism. The \emph{Procesi bundle} is defined to be 
$$\mathcal{P}_n = \rho_{n*}\OO_{X_n}$$
where $\rho_n:X_n\to H_n$ is the natural projection. 
Haiman's theorem \cite{Haiman} shows that $\rho_n$ is finite flat of degree $n!$. Thus, $\mathcal{P}_n$ is a vector bundle of rank $n!$. 

For a partition $\mu\vdash n$, its associated monomial ideal is 
$$I_\mu = \langle x^ay^b: (a,b)\notin \mu\rangle\in H_n. $$
Identify the scheme-theoretic fiber of $\rho_n$ at $I_\mu$ with a subscheme of $(\mathbb{A}^2)^n$, and denote its ideal by $J_\mu$. Then 
$$\mathbb{C}[x_1,y_1,\ldots,x_n,y_n]/J_\mu
\simeq \mathcal{P}_n(I_\mu). $$
By Haiman \cite[Proposition 3.7.3]{Haiman}, 
$$J_\mu^\perp = 
\{f\in \mathbb{C}[x_1,y_1,\ldots,x_n,y_n]:
\langle f,g\rangle =0
\text{ for all $g\in J_\mu$}
\}=\mathsf{GH}_\mu,$$
where the pairing on $\mathbb{C}[x_1,y_1,\ldots,x_n,y_n]$ is given by 
$$\langle f,g\rangle = f(
\tfrac{\partial}{\partial x_1},
\ldots, 
\tfrac{\partial}{\partial y_n})g(x_1,
\ldots, 
y_n)\big|_{x_1=
\cdots=y_n=0}.$$

The torus $T=(\mathbb{C}^\times)^2$ acts on $\mathbb{A}^2$ with weight $q$ on $x$, and weight $t$ on $y$.
The schemes, sheaves, and vector spaces above admit actions of the torus $T=(\mathbb{C}^\times)^2$ and the symmetric group $S_n$. The latter action is induced by the permutation of the points in $(\mathbb{A}^2)^n$. 
Moreover, Haiman \cite{Haiman} proved 
$$
\operatorname{Frob}_{q,t}(\mathcal{P}_n(I_\mu))= 
\operatorname{Frob}_{q,t}(\mathsf{GH}_\mu) 
= \widetilde{H}_\mu[X;q,t].$$
Here, for a finite-dimensional bigraded $S_n$-module $U$, 
\[\operatorname{Frob}_{q,t}(U):=\sum_{a,b}q^at^b\sum_{\mu \vdash n}\dim\Hom_{S_n}(V_\mu, U_{a,b})s_\mu(X),\]
where $V_\mu$ is the irreducible representation of $S_n$ indexed by $\mu$ and $s_\mu(X)$ is the Schur function.

\subsection{Nested ideals}
Fix $n\geq 2$. 
The \emph{nested Hilbert scheme} is 
$$H_{n-1,n}
=\{(I_1,I_2):I_1\subset I_2\}
\subset H_{n}\times H_{n-1}.$$
By a theorem of Tikhomirov (see \cite[Proposition 3.5.2]{Haiman} and \cite{Cheah}), $H_{n-1,n}$ is a smooth irreducible variety of dimension $2n$. 
Let $q_{n-1}:H_{n-1,n}\to H_{n-1}$ and $p_{n}:H_{n-1,n}\to H_n$ be the two projections. 
We have a tautological exact sequence on $H_{n-1,n}$
\begin{equation}
\label{eq:tau-N-Hn}
0
\longrightarrow
\mathcal{L} 
\longrightarrow
p_n^*\mathcal{B}_n 
\longrightarrow
q_{n-1}^*\mathcal{B}_{n-1}
\longrightarrow
0,
\end{equation}
where $\mathcal{L}$ is the line bundle whose fiber at $(I_1,I_2)$ is the quotient $I_2/I_1$. 

We have a morphism $\alpha: H_{n-1,n}\to U_n$ sending $(I_1,I_2)$ to $(I_1,a)$ where $\{a\} = \operatorname{Supp}(I_2/I_1)$. 
The projection $p_n:H_{n-1,n}\to H_n$ factors as $\pi_n\circ \alpha$. 

The fiber of $\alpha$ at $(I,a)\in U_n$ can be identified with 
$$\alpha^{-1}(I,a)=\mathbb{P}(\operatorname{soc}((\mathbb{C}[x,y]/I)_a)),$$
where the socle of a local ring $(R,\mathfrak{m})$ is 
$$\operatorname{soc}R := 
\{r\in R: \mathfrak{m}\cdot r =0\}
=\Hom_R(R/\mathfrak{m},R).$$
By Ellingsrud and Str\o mme \cite{EllStr}, this can be globalized as 
\begin{equation}
\label{eq:nHil-is-proj}
H_{n-1,n} = \mathbb{P}^\vee_{U_n}(\kappa_{U_n}):= 
\operatorname{Proj}_{U_n}
\big(\operatorname{Sym}(\kappa_{U_n})\big).
\end{equation}
Here $\kappa_{U_n}$ is the dualizing sheaf of $U_n$. 
Moreover, $\mathcal L$ in \eqref{eq:tau-N-Hn} is its tautological line bundle $\mathcal O(-1)$. 

\begin{remark}\label{rmk:relativeshv}
Let us briefly explain the geometric intuition here. 
Note that coherent sheaves on $U_n$ can be identified with coherent sheaves on $H_n$ equipped with a $\mathcal{B}_n$-module structure. 
The relative dualizing sheaf $\kappa_{\pi_n}$ of $\pi_n:U_n\to H_n$ corresponds to the $\mathcal{B}_n$-module 
$\HHom_{H_n}(\mathcal{B}_n,\OO_{H_n})$. 
Its fiber at $(I,a)\in U_n$ can be identified with 
$$(\mathbb{C}[x,y]/I)^\vee/
\mathfrak{m}_a(\mathbb{C}[x,y]/I)^\vee
=(\operatorname{soc}((\mathbb{C}[x,y]/I)_a))^\vee.$$
Note that $\mathbb{P}^\vee (\kappa_{U_n})=\mathbb{P}^\vee(\kappa_{\pi_n})$ since $\kappa_{U_n}$ and $\kappa_{\pi_n}$ differ by tensoring with the line bundle $\pi_n^* \kappa_{H_n}$. 
Finally, the fiber of $\mathbb{P}^\vee(\cdots)$ parametrizes one-dimensional quotients. 
\end{remark}

\subsection{Local structure of $\alpha$}\label{sec:locstrofalpha}
By \eqref{eq:nHil-is-proj} and the Nakayama lemma, $\kappa_{U_n}$ is generated by $s:=\dim_{\mathbb{C}} \operatorname{soc}(\mathbb{C}[x,y]/I)_a$ elements near $\xi = (I,a)\in U_n\subset H_n\times \mathbb{A}^2$. 
This defines a closed embedding $(H_{n-1,n})_{\xi}$ into $\mathbb{P}^{s-1}_{H_n\times \mathbb{A}^2,\xi}$. 
We obtain the following diagram:
$$
\xymatrix@!R=2pc{
\mathbb{P}^{s-1}\ar[r]^-{\subset}
\ar[d]&
(H_{n-1,n})_{\xi}
\ar[r]^-{\subset}_{(*)}
\ar[d]^{\alpha}
& 
\mathbb{P}^{s-1}_{H_n\times \mathbb{A}^2,\xi}\ar[d]\\
\xi\ar[r]^-{\in}&
(U_n)_{\xi}
\ar[r]^-{\subset}
&(H_n\times \mathbb{A}^2)_\xi
}
$$
Song \cite[Claim B]{Song} proved that $(*)$ can be made a regular embedding. 
Let us recall the argument. 

\begin{remark}
Let us first recall the Hilbert--Burch theorem \cite[Section 20.4]{Eis}. 
Let $R$ be a regular local ring. 
Assume the quotient $R/\mathfrak{a}$ is Cohen--Macaulay of codimension $2$. By the Auslander--Buchsbaum formula \cite[\href{https://stacks.math.columbia.edu/tag/090V}{090V}]{Stacks}, we have a minimal free resolution 
$$0\longrightarrow 
R^{\oplus b_1}
\stackrel{\Phi}\longrightarrow 
R^{\oplus b_2}
\stackrel{\Psi}\longrightarrow 
R \longrightarrow R/\mathfrak{a}
\longrightarrow 0.
$$
Minimality implies that the entries of $\Phi,\Psi$ are in $\mathfrak{m}$. 
Taking ranks over the fraction field $\operatorname{Frac}(R)$ gives $b_2=b_1+1$; see 
\cite[Corollary 20.13]{Eis}. 
Hilbert--Burch \cite[Theorem 20.15]{Eis} implies 
$$
(\text{the $i$-th entry of $\Psi$})
= 
\pm r
(\text{the $i$-th ($b_1\times b_1$)-minor of $\Phi$})$$
for some nonzerodivisor $r\in R$. 
If $r$ were not a unit, then $\mathfrak{a}\subset (r)$ would imply that the height of $\mathfrak{a}$ is at most $1$. 
Thus $r$ is a unit, and a change of basis allows us to assume $r=1$. 
\end{remark}

Now let $\xi=(I,a)\in U_n$. 
Let 
\begin{equation}\label{eq:Sam}
S = 
\OO_{H_n\times \mathbb{A}^2,\xi},\qquad 
\mathfrak{a}\trianglelefteq S\text{ with }
S/\mathfrak{a}= \OO_{U_n,\xi},\qquad 
S/\mathfrak{m}\simeq \mathbb{C}.
\end{equation}
Since $\pi_n$ is finite flat over the smooth scheme $H_n$, the universal family $U_n$ is also Cohen--Macaulay of dimension $2n$, i.e., codimension $2$ in $H_n\times \mathbb{A}^2$. 
So we can apply Hilbert--Burch to $S$ and $S/\mathfrak{a}$. 
Choose a minimal free resolution
$$0\longrightarrow 
S^{\oplus s}
\stackrel{\Phi}\longrightarrow 
S^{\oplus s+1}
\stackrel{\Psi}\longrightarrow 
S
\longrightarrow 
S/\mathfrak{a}\longrightarrow 0.
$$
Recall that the dualizing sheaf is given by
$$\kappa_{S/\mathfrak{a}} =\Ext_S^2(S/\mathfrak{a},S)
=\operatorname{coker}
[S^{\oplus s+1}
\stackrel{\Phi^\top}\longrightarrow S^{\oplus s}]. $$
Note that
$\operatorname{Sym}_S (S^{\oplus s})
\simeq S[T_1,\ldots,T_s]$, 
and we can write 
$$\operatorname{Sym}_S(\kappa_{S/\mathfrak{a}})
= S[T_1,\ldots,T_s]/\langle\eta_1,\ldots,\eta_{s+1} \rangle$$
where
$\eta_i = \sum_{j=1}^s \Phi_{ij}T_j$ for $1\leq i\leq s+1$. 
This realizes 
$$(H_{n-1,n})_{\xi} = \mathbb{P}^\vee_{S/\mathfrak{a}}(\kappa_{S/\mathfrak{a}})
\subset 
\mathbb{P}^\vee_S(S^{\oplus s}) 
=
\mathbb{P}^{s-1}_{H_n\times \mathbb{A}^2,\xi}
$$
as the zero locus of $\eta_1,\ldots,\eta_{s+1}\in H^0(\mathbb{P}^{s-1}_S,\OO(1))$. 
Since the entries of $\Phi$ are all in $\mathfrak{m}$, 
the fiber of $(H_{n-1,n})_\xi$ at $R/\mathfrak{m}$, isomorphic to $\mathbb{P}(\operatorname{soc}(\mathbb{C}[x,y]/I)_a)$, is exactly $\mathbb{P}^{s-1}$. 
So we have 
$$s=\dim \operatorname{soc}(\mathbb{C}[x,y]/I)_a.$$
The scheme $(H_{n-1,n})_{\xi}$ has dimension $2n$ and codimension $s+1$ in $\mathbb{P}^{s-1}_S$.

Hence $\eta_1,\ldots,\eta_{s+1}$ form a regular sequence. 
We thus obtain the Koszul resolution on $\mathbb{P}_S^{s-1}$:
\begin{equation}
\label{eq:KoszulfornHil}
\begin{aligned}
0 & 
\longrightarrow 
\OO_{\mathbb{P}^{s-1}_S}(-s-1)
\longrightarrow 
\OO_{\mathbb{P}^{s-1}_S}(-s)
^{\oplus\binom{s+1}{s}}
\longrightarrow \cdots \\
& \longrightarrow 
\OO_{\mathbb{P}^{s-1}_S}(-1)^{\oplus \binom{s+1}{1}}
\longrightarrow 
\OO_{\mathbb{P}^{s-1}_S}
\longrightarrow 
\OO_{(H_{n-1,n})_{\xi}}
\longrightarrow 0. 
\end{aligned}
\end{equation}
This will be useful in the next section.

\section{A pushforward theorem for mixed polygraphs}\label{sec:polygraph}

For any $b\in \mathbb{Z}_{\geq 0}$, we construct the following diagram:
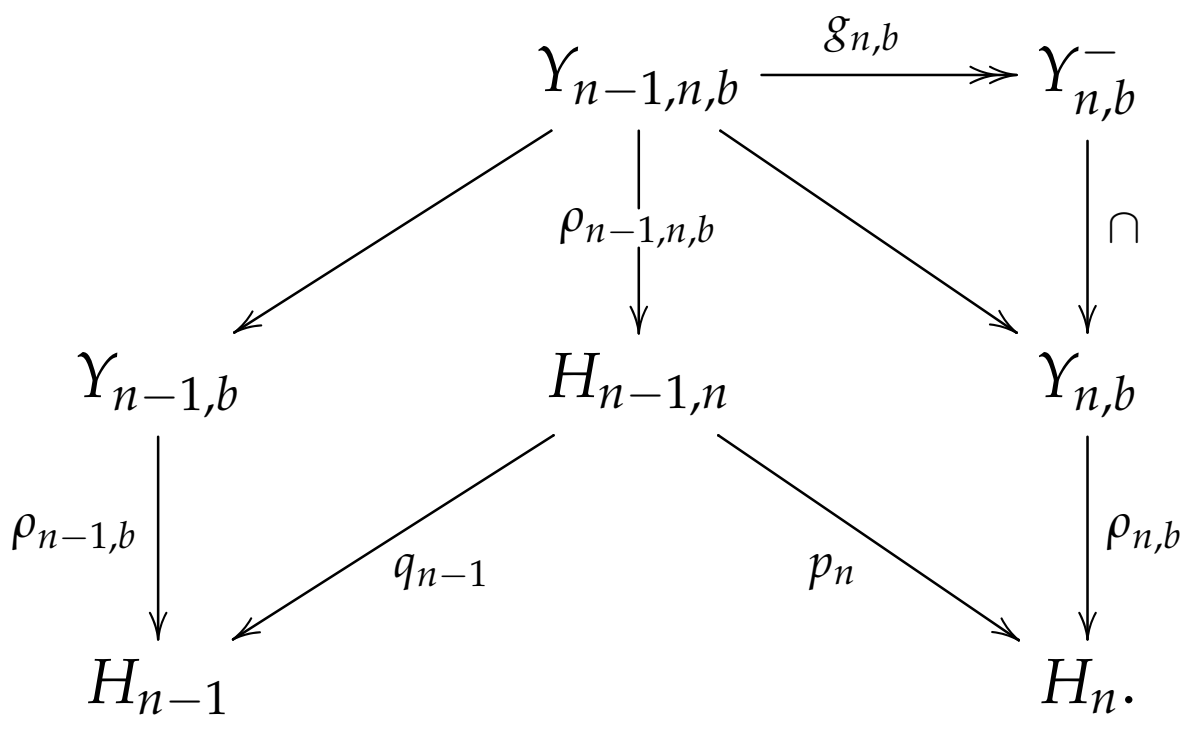
\begin{equation}
\label{eq:mixed-main-diagram}
\begin{matrix}
\xymatrix@!R=2.8pc@C=4pc{
& Y_{n-1,n,b}
    \ar[d]|{\rho_{n-1,n,b}}
    \ar@{->>}[r]^{g_{n,b}}
    \ar[dl]
    \ar[dr] & 
Y_{n,b}^-
    \ar[d]^{\cap}\\
Y_{n-1,b}
    \ar[d]_{\rho_{n-1,b}}& 
H_{n-1,n}
    \ar[dl]^{q_{n-1}}
    \ar[dr]_{p_n}&
Y_{n,b}
    \ar[d]^{\rho_{n,b}}\\
H_{n-1}&&H_n.}
\end{matrix}
\end{equation}
Let 
$$Y_{n,b} = X_n\times_{H_n} 
\underbrace{
U_n
\times_{H_n}\cdots\times_{H_n}U_n}_{\text{$b$ factors}}$$
and let $\rho_{n,b}:Y_{n,b}\to H_n$ be the natural projection.
The map $\rho_{n,b}$ is finite flat of degree $n!\cdot n^b$ with 
$$\rho_{n,b,*}\OO_{Y_{n,b}} = \mathcal{P}_n\otimes \mathcal{B}_n^{\otimes b}.$$
Since 
$Y_{n,b}\to X_n$ is finite flat of degree $n^b$, with reduced generic fiber, for any $p\in X_n$, 
$$\OO_{Y_{n,b},p}
\hookrightarrow 
\OO_{Y_{n,b},p}\otimes_{\OO_{X_n,p}}
\operatorname{Frac}(\OO_{X_n,p})
=
\operatorname{Frac}(\OO_{X_n,p})^{\Pi n^b}.
$$
So $Y_{n,b}$ is reduced, and we have the following explicit description:
\begin{align*}
Y_{n,b}& = 
\left\{(I,P_1,\ldots,P_n,Q_1,\ldots,Q_b):
\text{each }Q_i\in \{P_1,\ldots,P_n\}
\right\}\\
& = \bigcup_{f:[b]\to [n]}
\left\{(I,P_1,\ldots,P_n,Q_1,\ldots,Q_b):
Q_i=P_{f(i)}
\right\}
\subset X_n\times (\mathbb{A}^2)^b.
\end{align*}
Note that this is a decomposition into irreducible components, and each component is isomorphic to $X_n$ via the natural projection. 

Next, let us define a reduced subscheme of $Y_{n,b}$:
\begin{align*}
Y_{n,b}^-& = 
\left\{(I,P_1,\ldots,P_n,Q_1,\ldots,Q_b):
\text{each } Q_i
\in \{P_1,\ldots,P_{n-1}\}
\right\}\\
& = \bigcup_{f:[b]\to [n-1]}
\left\{(I,P_1,\ldots,P_n,Q_1,\ldots,Q_b):
\begin{matrix}
Q_i=P_{f(i)}
\end{matrix}
\right\}
\subset Y_{n,b}.
\end{align*}

Let us define 
$$Y_{n-1,n,b} = H_{n-1,n}\times_{H_{n-1}} Y_{n-1,b}$$
with the natural morphism $\rho_{n-1,n,b}:Y_{n-1,n,b}\to H_{n-1,n}$. 
Since $\rho_{n-1,b}:Y_{n-1,b}\to H_{n-1}$ is finite flat, with reduced generic fiber, $Y_{n-1,n,b}$ is reduced, and 
it is the variety of tuples
\begin{equation}\label{eq:reprofYnnb}
(I_1\subset I_2,P_1,\ldots,P_{n-1},Q_1,\ldots,Q_b)
\end{equation}
such that 
$(I_1\subset I_2)\in H_{n-1,n}$,
$(I_2,P_1,\ldots,P_{n-1})\in X_{n-1}$,
and each $Q_i \in \{P_1,\ldots,P_{n-1}\}$.
The natural morphism $\alpha:H_{n-1,n}\to U_n$ induces a surjective map
$$g_{n,b}:Y_{n-1,n,b}\to Y_{n,b}^-$$ by sending \eqref{eq:reprofYnnb} to 
\begin{equation}\label{eq:imofYnnb}
(I_1,P_1,\ldots,P_{n-1},a,Q_1,\ldots,Q_b),
\end{equation}
where $\{a\} = \operatorname{Supp}(I_2/I_1)$. 
We have an isomorphism
\begin{equation}
\label{eq:Ynnbasfiber}
Y_{n-1,n,b}
\stackrel{\sim}\longrightarrow
\big(H_{n-1,n}\times_{U_n} Y_{n,b}^-\big)_{\rm red}
\hookrightarrow 
H_{n-1,n}\times_{U_n} Y_{n,b}^-,
\end{equation}
where $Y_{n,b}^-\to U_n$ is the map sending $(I,P_1,\ldots,P_n,Q_1,\ldots,Q_b)$ to $(I,P_n)$.
Since the fiber of $\alpha$ is reduced, 
\begin{equation}\label{eq:fiberofgnb}
g_{n,b}^{-1}(\eqref{eq:imofYnnb})
=\alpha^{-1}(I_1,a).
\end{equation}

The purpose of this section is to show the following. 

\begin{maintheorem}\label{thm:Section3}
The natural map induces
$$
\OO_{Y_{n,b}^-}
\stackrel{\sim}\longrightarrow 
Rg_{n,b,*}\OO_{Y_{n-1,n,b}}.
$$
\end{maintheorem}

Haiman \cite{Haiman} proved the case $b=0$, where $Y_{n,0}^-=Y_{n,0}=X_n$ and $Y_{n-1,n,0}$ was denoted by $X_{n-1,n}$ \textit{loc. cit.}. 
The approach in this section is based on his proof. 

\subsection{Haiman's criterion and induction}\label{sec:pushforward}
The following criterion is from Haiman \cite[Lemma 3.8.5]{Haiman}, with its fiber-dimension bound replaced by the vanishing it supplies.

\begin{theorem}[Haiman]
\label{thm:haiman-pushforward}
Let $f:Y\to X$ be proper between finite-type complex schemes.
For $z_1,\ldots,z_c\in H^0(X,\OO_X)$, let
$Z$ be the zero locus of $z_1,\ldots,z_c$ and $U=X\setminus Z$. 
Assume:
\begin{enumerate}
\item \label{enm:h-cri-i}
The sequence $z_1,\ldots,z_c$ is regular along $Z$;
\item \label{enm:h-cri-ii}
The sequence $f^*z_1,\ldots,f^*z_c$ is regular along $f^{-1}(Z)$;
\item \label{enm:h-cri-iii}
$R^jf_*\OO_Y=0$ for $j\geq c-1$;
\item \label{enm:h-cri-iv}
$\OO_U\to R(f|_{f^{-1}(U)})_*\OO_{f^{-1}(U)}$ is an isomorphism.
\end{enumerate}
Then $\OO_X\stackrel{\sim}\longrightarrow Rf_*\OO_Y$ by the natural morphism.
\end{theorem}
\begin{proof}
Without loss of generality, we may assume that $X$ is affine. 
For $(X,Z)$ or
$(Y,f^{-1}(Z))$, the regularity assumptions \eqref{enm:h-cri-i} and \eqref{enm:h-cri-ii} give vanishing of the local cohomology
sheaves in degrees $j<c$. The local-to-global spectral sequence
then gives 
$$
H_{Z}^j(X,\OO_X)=0,\qquad 
H_{f^{-1}(Z)}^j(Y,\OO_Y)=0,\qquad 
(j<c). $$
For $0\leq j<c-1$, the restriction sequence gives 
\[
 \underbrace{H_Z^j(X,\OO_X)}_{=0}
 \longrightarrow
 H^j(X,\OO_X)
 \stackrel{\sim}\longrightarrow
 H^j(U,\OO_X)
 \longrightarrow
 \underbrace{H_Z^{j+1}(X,\OO_X)}_{=0};
\]
and similarly, $H^j(Y,\OO_Y)
\simeq H^j(f^{-1}(U),\OO_Y)$. 
Now \eqref{enm:h-cri-iv} gives, for $0\leq j<c-1$,
\[
 H^j(Y,\OO_Y)\simeq H^j(f^{-1}(U),\OO_Y)
 \simeq H^j(U,\OO_U)\simeq H^j(X,\OO_X).
\]
For $j\geq c-1$, \eqref{enm:h-cri-iii} gives
$H^j(Y,\OO_Y)=\Gamma(X,R^jf_*\OO_Y)=0$.
\end{proof}

We will apply Haiman's criterion \ref{thm:haiman-pushforward} to 
$f=g_{n,b}$ with 
\begin{equation}\label{eq:applyHc}
z_i = y_i-y_{i+1},\qquad 1\leq i\leq c=n-1,
\end{equation}
where $y_i$ is the second coordinate of $P_i$ for $(I,P_1,\ldots,P_n,Q_1,\ldots,Q_b)\in Y_{n,b}^-$. 

\subsection{Condition \eqref{enm:h-cri-i}}
A map $f:[b]\to [n]$ induces a morphism 
$f: (\mathbb{A}^2)^n\longrightarrow (\mathbb{A}^2)^b$. 
Let $\Gamma_f\subset (\mathbb{A}^2)^{b+n}$ be its graph.
Let 
$$\Gamma_{n,b}^- = \bigcup_{f:[b]\to [n-1]} \Gamma_f
=\{(P_1,\ldots,P_n,Q_1,\ldots,Q_b):
Q_i\in \{P_1,\ldots,P_{n-1}\}\}.$$
Let 
$$A(n,b)^-=
\frac{\mathbb{C}[x_1,y_1,\ldots,x_n,y_n,u_1,v_1,\ldots,u_b,v_b]}{
\bigcap_{f:[b]\to [n-1]}
\langle u_i-x_{f(i)},v_i-y_{f(i)}:1\leq i\leq b\rangle}
$$
be the coordinate ring of $\Gamma_{n,b}^-$. 
Let 
$\mathfrak{alt}_n
\trianglelefteq 
\mathbb{C}[x_1,\ldots,y_n]$ be the ideal generated by alternating polynomials (i.e., polynomials $f$ such that $\sigma f=(-1)^{\ell(\sigma)}f$ for $\sigma \in S_n$). 
Haiman \cite[Proposition 3.4.2]{Haiman} proved that 
\begin{equation}\label{eq:Haimanblowup}
X_n = \operatorname{Bl}_{\mathfrak{alt}_n}
(\mathbb{A}^2)^n=
\operatorname{Proj}\left(
\bigoplus_{m\geq 0}\mathfrak{alt}_n^m \cdot t^m\right).
\end{equation}

\begin{proposition}
\label{prop:Ynbisblowup}
We can realize $Y_{n,b}^-$ as a blow-up
\begin{equation}\label{eq:myblowup}
Y_{n,b}^- = \operatorname{Bl}_{\mathfrak{alt}_nA(n,b)^-}
\Gamma_{n,b}^-
= \operatorname{Proj}\left(
\bigoplus_{m\geq 0}\mathfrak{alt}_n^m A(n,b)^-\cdot t^m\right).
\end{equation}
\end{proposition}
\begin{proof}
Let $B$ be the blow-up in \eqref{eq:myblowup}. 
Since we have a natural surjective map
$$
A(n,b)^-\otimes_{\mathbb{C}[x_1,\ldots,y_n]}
\left(
\bigoplus_{m\geq 0}\mathfrak{alt}_n^m \cdot t^m\right)
\longrightarrow 
\bigoplus_{m\geq 0}\mathfrak{alt}_n^m A(n,b)^-\cdot t^m,
$$
we have a closed embedding 
$B\stackrel{\subset}\longrightarrow X_n\times_{(\mathbb{A}^2)^n} \Gamma_{n,b}^-$.
Since $A(n,b)^-[t]$ is reduced, $B$ is a reduced scheme. 
Note that $$(X_n\times_{(\mathbb{A}^2)^n} \Gamma_{n,b}^-)_{\rm red}\simeq Y_{n,b}^-.$$
It suffices to show $B$ contains all generic points of irreducible components of $(X_n\times_{(\mathbb{A}^2)^n} \Gamma_{n,b}^-)_{\rm red}$. 
Restrict to the open locus $U\subset (\mathbb{A}^2)^n$ where $y_i\neq y_j$ for $i\neq j$. 
Now the unions are all disjoint: 
$$
(\Gamma_{n,b}^-)_U = \coprod_{f:[b]\to [n-1]} 
(\Gamma_{f})_U,\qquad B_U = \coprod_{f:[b]\to [n-1]} \operatorname{Bl}_{\mathfrak{alt}_n}
(\Gamma_f)_U,$$
with $(\Gamma_f)_U\simeq U$.
Since
$$\operatorname{Bl}_{\mathfrak{alt}_n}
(\Gamma_f)_U
=\left\{(I,P_1,\ldots,P_n,Q_1,\ldots,Q_b):
Q_i=P_{f(i)}
\right\}_U,$$
the assertion follows. 
\end{proof}

\begin{remark}
Let 
$\Gamma_{n,b} = \bigcup_{f:[b]\to [n]} \Gamma_f$ and define $A(n,b)$ as above. We can prove 
$$Y_{n,b} = \operatorname{Bl}_{\mathfrak{alt}_nA(n,b)}\Gamma_{n,b}=
\operatorname{Proj}\left(
\bigoplus_{m\geq 0}\mathfrak{alt}_n^m A(n,b)\cdot t^m\right).$$
\end{remark}

The proof of the following proposition generalizes the argument of Haiman \cite[Proposition 4.2]{HaimanNotes}. 

\begin{proposition}
\label{prop:altA-isfree}
The ideal $\mathfrak{alt}_n^m A(n,b)^-$ is a free $\mathbb{C}[y_1,\ldots,y_n]$-module. 
\end{proposition}

To prove this proposition, we need Haiman's common-basis theorem \cite[Theorem 4.2]{Haiman}. 

\begin{theorem}
\label{thm:haiman-common-basis}
Let $0\leq r\leq n$, $0\leq k\leq b$ and $a\in \mathbb{Z}$. 
Set
$$
Y_{n,b}(a,r,k)=
\bigcup_{\substack{f:[b]\to[n],\ T\subseteq[n]\\
|(T\cap[r])\setminus f([k])|\geq a}}
\big(\{x_j=0\}_{j\in T}\cap 
\Gamma_f \big)
\subseteq \Gamma_{n,b}. $$
Let
$I_{n,b}(a,r,k)\trianglelefteq A(n,b)$ 
be its ideal. 
There is a $\mathbb{C}[y_1,\ldots,y_n]$-basis $\mathfrak{B}$ of $A(n,b)$ such that $\mathfrak{B}\cap I_{n,b}(a,r,k)$ forms a basis of $I_{n,b}(a,r,k)$. 
In particular, the coordinate ring
$$A(n,b)/I_{n,b}(a,r,k)$$
of $Y_{n,b}(a,r,k)$ is free over $\mathbb{C}[y_1,\ldots,y_n]$.
\end{theorem}

Informally, we can describe 
$Y_{n,b}(a,r,k)=
\bigcup_{f:[b]\to [n]}
Y_{n,b}(a,r,k)_f$, where $Y_{n,b}(a,r,k)_f\subseteq \Gamma_f$ is the locus on which
at least $a$ elements (counting multiplicity) of
$\{x_1,\ldots,x_r\}\setminus 
\{u_{1},\ldots,u_{k}\}$ vanish.

\begin{proof}[Proof of Proposition \ref{prop:altA-isfree}]
Put $B=b+nm$. 
Let us consider 
$$L = 
\bigcup_{f\in F}
\Gamma_f,\qquad 
F = \{
[B]\stackrel{f}\to [n]:
f([b])\subseteq [n-1]\}.$$
Note that 
$$f\in F\iff 
|(T\cap \{n\})\setminus f([b])|\geq 1
\text{ for some }T\subseteq [n],$$
and it happens only when $n\in T$.
By relabeling, we have 
$$L \cap \{x_n=0\} \simeq Y_{n,B}(1,1,b). $$
By translating the first coordinate, we get 
$$L\simeq \mathbb{A}^1\times Y_{n,B}(1,1,b). $$
By Theorem \ref{thm:haiman-common-basis}, $\mathbb{C}[L]$ is free over $\mathbb{C}[y_1,\ldots,y_n]$. 

Write the points of $L$ as 
$$(P_1,\ldots,P_n,Q_1,\ldots,Q_b,
Q_{11},\ldots,Q_{1n},\dots 
Q_{m1},\ldots,Q_{mn}).$$
Let $G_m = S_n^m\subset S_{nm}$ act on $(\mathbb{A}^{2})^{B}=(\mathbb{A}^{2})^{b+nm}$ by permuting the last $nm$ coordinates $Q_{ij}$.  
Let $e\in \mathbb{C}[G_m]$ be the idempotent for the tensor product of the sign representations of $G_m$. 
We will prove 
$$
e\cdot \mathbb{C}[L]
\stackrel{\sim}\longrightarrow \mathfrak{alt}_n^m A(n,b)^-.$$
Since $\mathbb{C}[L]$ is free over $\mathbb{C}[y_1,\ldots,y_n]$, so is $e\cdot \mathbb{C}[L]$; thus Proposition \ref{prop:altA-isfree} follows. 

We have a morphism 
$\epsilon:\Gamma_{n,b}^-\to L$ 
sending $(P_1,\ldots,P_n,Q_1,\ldots,Q_b)$ to 
$$(P_1,\ldots,P_n,Q_1,\ldots,Q_b,
P_1,\ldots,P_n,
\stackrel{m}{\ldots\ldots},
P_1,\ldots,P_n).$$
That is, $Q_{ij}=P_j$. 
This induces $\epsilon^*:\mathbb{C}[L]\to A(n,b)^-$. 

First, we show that $\epsilon^*(e\cdot \mathbb{C}[L])=\mathfrak{alt}_n^m A(n,b)^-$.
Note that $e\cdot \mathbb{C}[L]$ is generated by 
$$\varphi_0(P_1,\ldots,Q_b)\cdot 
\varphi_1(Q_{11},\ldots,Q_{1n})\cdots
\varphi_m(Q_{m1},\ldots,Q_{mn}),
$$
with $\varphi_0\in A(n,b)^-$ and  $\varphi_1,\ldots,\varphi_m$ alternating. 
Its image under $\epsilon^*$ is 
$$\varphi_0(P_1,\ldots,Q_b)\cdot 
\varphi_1(P_{1},\ldots,P_{n})\cdots 
\varphi_m(P_{1},\ldots,P_{n})
\in \mathfrak{alt}_n^mA(n,b)^-. 
$$
These products generate $\mathfrak{alt}_n^mA(n,b)^-$. 
Let us prove the induced map is injective. Assume $\varphi\in e\cdot \mathbb{C}[L]$. 
Note that $\varphi|_{\Gamma_f}=0$ unless
$$f(b+(k-1)n+1),\ldots,f(b+kn)$$ 
are distinct for all $1\leq k\leq m$. 
The union of the $\Gamma_f$ with this property is exactly the $G_m$-orbit of the image of $\epsilon$. 
Since $\varphi$ is $G_m$-alternating, $\varphi$ is determined by its restriction to the image of $\epsilon$. 
This proves the injectivity. 
\end{proof}

Propositions \ref{prop:Ynbisblowup} and \ref{prop:altA-isfree} give condition \eqref{enm:h-cri-i} of Haiman's criterion. 

\begin{corollary}
\label{cor:mixed-target-regular}
The scheme $Y^-_{n,b}$ is flat over $\mathbb{C}[y_1,\ldots,y_n]$.
In particular, the sequence
$z_1,\ldots,z_{n-1}$ is regular.
\end{corollary}
\begin{proof}
Since flatness implies torsion-freeness, each $z_i$ is a nonzerodivisor in the successive quotient. 
\end{proof}

\subsection{Condition \eqref{enm:h-cri-ii}}
The proof of condition \eqref{enm:h-cri-ii} is easier. 
Since $\rho_{n-1,n,b}:Y_{n-1,n,b}\to H_{n-1,n}$ is finite flat, 
$Y_{n-1,n,b}$ is Cohen--Macaulay. 
It suffices to examine the dimension.

\begin{proposition}
\label{prop:mixed-source-regular}
The sequence $g_{n,b}^*z_1,\ldots,g_{n,b}^*z_{n-1}$ is regular. 
\end{proposition}
\begin{proof}
Let $X_{n-1,n}:=Y_{n-1,n,0}$ and 
$r:Y_{n-1,n,b}\longrightarrow X_{n-1,n}$ be the morphism 
forgetting the $b$ attached points. 
By definition, it is finite flat, so it suffices to treat the case $b=0$. 
This is proved by Haiman \cite[Proposition 3.5.6]{Haiman}. 
\end{proof}

\subsection{Condition \eqref{enm:h-cri-iii}}
Since Theorem \ref{thm:haiman-pushforward} is stated with its fiber-dimension bound replaced by the vanishing it supplies, 
Haiman's argument simplifies slightly: 
we no longer need the explicit computation for $n=3$. 

\begin{proposition}
\label{cor:mixed-cohomological-bound}
For any $\eta\in Y_{n,b}^-$, 
we have a closed embedding
$$(Y_{n-1,n,b})_\eta
\subseteq \mathbb{P}^{n-2}_{(Y_{n,b}^-)_{\eta}}.$$
In particular, when $n=2$, $g_{n,b}$ is an isomorphism, and when $n>2$, 
$$R^jg_{n,b,*}\mathcal{O}_{Y_{n-1,n,b}} = 0 \text{ for $j\geq n-2>0$}.$$
\end{proposition}
\begin{proof}
By \eqref{eq:Ynnbasfiber} and \eqref{eq:nHil-is-proj}, we can embed
$$Y_{n-1,n,b} 
\stackrel{\subset}\longrightarrow 
H_{n-1,n}\times_{U_n} Y_{n,b}^-
=
\mathbb{P}_{Y_{n,b}^-}(\kappa_{U_n}).$$
By Nakayama's lemma, $\kappa_{U_n}$ is generated by $s:=\dim_{\mathbb{C}} \operatorname{soc}(\mathbb{C}[x,y]/I)_a$ elements near $\eta=:(I,P_1,\ldots,a,Q_1,\ldots,Q_b)$. 
This defines a surjective morphism 
$\mathcal{O}_{Y_{n,b}^-,\eta}^{\oplus s}\longrightarrow \kappa_{U_n,\eta}$, 
which induces
$$(Y_{n-1,n,b})_\eta
\subseteq \mathbb{P}^{s-1}_{(Y_{n,b}^-)_{\eta}}.$$
So the assertion follows from 
$s\leq n-1$ as $n\geq 2$. 

Over $X=\mathbb{P}^{n-2}_{(Y_{n,b}^-)_{\eta}}$, 
we have a surjective map 
$\mathcal{O}_{X}\twoheadrightarrow\mathcal{O}_{(Y_{n-1,n,b})_\eta}.$
Since $H^{n-2}(X,\mathcal{O}_X)=0$, we have 
$H^{n-2}((Y_{n-1,n,b})_\eta,\mathcal{O}_{(Y_{n-1,n,b})_\eta})=0$. 
\end{proof}

\subsection{Condition \eqref{enm:h-cri-iv}}
The last condition follows from induction.

For any partition $[n]=I\sqcup J$ with $n\in I$, consider the open subset
$$\Omega_{I,J} = 
 \{P_i\neq P_j\text{ for }i\in I,j\in J\}
 \subseteq Y_{n,b}^-.$$
Note that 
$$\bigcup_{\substack{[n]=I\sqcup J \\|J|\geq 1}} \Omega_{I,J}
= \text{complement of }\{P_1=\cdots =P_n\}$$
covers the complement of the zero locus of $z_1,\ldots,z_{n-1}$.

\begin{proposition}
\label{prop:mixed-open-comparison}
We have 
$$\xymatrix@C=0pc{
&(Y_{n-1,n,b})_{\Omega_{I,J}}
\ar[rr]^-{g_{n,b}} 
\ar@{<-}[d]^{\wr}
&&
(Y_{n,b}^-)_{\Omega_{I,J}}
\ar@{<-}[d]^{\wr}\\
\displaystyle
\coprod_{[b]=L\sqcup K}
\bigg(&
\big(Y_{|I|-1,|I|,|L|}\times 
Y_{|J|,|K|}\big)_{\rm disj}
\ar[rr]^{
g_{|I|,|L|}
\times \operatorname{id}}
&\qquad\qquad\qquad&
\big(Y_{|I|,|L|}^{-}\times 
Y_{|J|,|K|}\big)_{\rm disj}&\bigg).
}$$
The subscript $\mathrm{disj}$ denotes the open subset on which the supports of the two finite subschemes are disjoint.
\end{proposition}
\begin{proof}
Note that $(Y_{n-1,n,b})_{\Omega_{I,J}}$ and 
$(Y_{n,b})_{\Omega_{I,J}}$ can be written as a disjoint union of 
$$\coprod_{[b]=L\sqcup K}\bigg(
\bigcup_{
\substack{
f:[b]\to [n]\text{ with }\\
f:L\to I, f:K\to J}
}\{\cdots\}
\bigg).$$
The left vertical morphism sends 
$$
\big(I_1\subset I_2,(P_i)_{i\in I\setminus\{n\}},(Q_\ell)_{\ell\in L}\big)
\times 
\big(I_3,(P_j)_{j\in J},(Q_k)_{k\in K}\big)$$
to 
$$(I_1I_3\subset I_2I_3, P_{1},\ldots,P_{n-1},Q_1,\ldots,Q_b).$$
Similarly for the right vertical map. 
Commutativity is left to the reader.
\end{proof}

\subsection{Proof of Theorem \ref{thm:Section3}}

Induct on $n\geq 2$, for all $b\geq0$ simultaneously.
Proposition \ref{cor:mixed-cohomological-bound} proves the case $n=2$.
For $n\geq3$, apply Haiman's criterion \ref{thm:haiman-pushforward} with \eqref{eq:applyHc}. 
Its four conditions are
\[
\begin{array}{ll}
 \text{\eqref{enm:h-cri-i}}&\text{Corollary \ref{cor:mixed-target-regular}},\\
 \text{\eqref{enm:h-cri-ii}}&\text{Proposition \ref{prop:mixed-source-regular}},\\
 \text{\eqref{enm:h-cri-iii}}&\text{Proposition \ref{cor:mixed-cohomological-bound}},\\
 \text{\eqref{enm:h-cri-iv}}&\text{Proposition \ref{prop:mixed-open-comparison}
                      and the induction hypothesis}.
\end{array}
\]
The criterion therefore proves the isomorphism.

\section{The Procesi bundle near a socle fiber}\label{sec:tilting}

Let $\xi=(I,a)\in U_n$.
Recall the construction in Section \ref{sec:locstrofalpha}. 
In this section, let us denote 
\begin{equation}\label{eq:XYZ}
\mathbb{P}^{s-1}\simeq Z=\alpha^{-1}(\xi)\subset 
Y := (H_{n-1,n})_\xi\subset 
X := \mathbb{P}^{s-1}_{H_n\times \mathbb{A}^2,\xi}.
\end{equation}
Let $\mathcal{Q}_X$ (resp. $\mathcal{Q}_Z$) be the universal quotient bundle over $X=\mathbb{P}^{s-1}_{H_n\times \mathbb{A}^2,\xi}$ (resp. $Z\simeq \mathbb{P}^{s-1}$). 
Define the following sheaves on 
$X$, $Y$ and $Z$:
\begin{equation}\label{eq:tilting}
\mathcal{T}_X=
{\textstyle \bigwedge^\bullet \mathcal{Q}_X}
:=\bigoplus_{i=0}^{s-1} {\textstyle \bigwedge^i \mathcal{Q}_X}
\qquad 
\mathcal{T}_Y = \mathcal{T}_X|_{Y},\qquad
\mathcal{T}_Z = \mathcal{T}_X|_{Z}
= {\textstyle \bigwedge^\bullet \mathcal{Q}_Z}.
\end{equation}

For a vector bundle $\mathcal T$, let $\operatorname{add}(\mathcal T)$ denote the bundles isomorphic to direct summands of finite direct sums of $\mathcal T$.

In this section, we will prove the following. 
\begin{maintheorem}\label{thm:Section4}
The vector bundle $q_{n-1}^*\mathcal{P}_{n-1}|_Y$ belongs to $\operatorname{add}(
\mathcal{T}_Y)$. 
In particular, $q_{n-1}^*\mathcal{P}_{n-1}|_Z$ belongs to $\operatorname{add}(
\mathcal{T}_Z)$. 
\end{maintheorem}

\subsection{Tilting bundles}
Let us first prove the following theorem. 
\begin{theorem}\label{thm:Tistilting}
The sheaf $\mathcal{T}_Y$ is a tilting generator for $D_{qc}(Y)$.  
\end{theorem}

\begin{remark}
Let us recall some tilting theory. 
For a quasi-compact, quasi-separated scheme $\mathfrak{X}$, let $D_{qc}(\mathfrak{X})$ be the unbounded derived
category of quasi-coherent sheaves. 
Let $\mathcal{T}$ be an object of $D_{qc}(\mathfrak{X})$. 
\begin{itemize}
    \item 
We say $\mathcal{T}$ is \emph{compact} if $\Hom_{D_{qc}(\mathfrak X)}(\mathcal T,-)$ commutes with arbitrary direct sums; see
\cite[\href{https://stacks.math.columbia.edu/tag/07LS}{07LS}]{Stacks}; 
by \cite[\href{https://stacks.math.columbia.edu/tag/09M1}{09M1}]{Stacks}, any (finite-rank) vector bundle over $\mathfrak{X}$ is compact. 
    \item 
We say $\mathcal{T}$ is a generator if, for every $\mathcal K\in D_{qc}(\mathfrak{X})$, 
$$R\Hom(\mathcal{T},\mathcal K)=0\Longrightarrow\mathcal K=0.$$
Equivalently, the smallest localizing triangulated subcategory containing $\mathcal T$, that is, the one closed under arbitrary direct sums, is $D_{qc}(\mathfrak X)$. 
    \item 
We say $\mathcal{T}$ is tilting if it is a compact generator and 
$$\Hom_{D_{qc}(\mathfrak{X})}(\mathcal{T},\mathcal{T}[j])=0
\text{ whenever }j\neq 0.$$
When $\mathcal{T}$ is a vector bundle that is already a generator, the remaining self-extension condition is equivalent to 
$$H^{>0}(\mathfrak{X},\EEnd_{\mathfrak{X}}(\mathcal{T}))=0.$$

\end{itemize}
\end{remark}

\begin{lemma}\label{lem:Lambdaistitling}
Over $X$, the vector bundle 
$\mathcal{T}_X$ is a compact generator of $D_{qc}(X)$, and we have 
$$H^k(X,\EEnd_X(\mathcal{T}_X)(-j))=0
\qquad 
\text{for any $0\leq j<k$}.$$
\end{lemma}
\begin{proof}
Recall the Euler exact sequence on $X$:
$$0
\longrightarrow 
\OO_X(-1)
\longrightarrow 
\OO_X^{\oplus s}
\longrightarrow 
\mathcal{Q}_X
\longrightarrow 
0.
$$
Taking exterior powers, we get 
an exact sequence
$$
\begin{aligned}
0&\longrightarrow \OO_X(-c)
\longrightarrow \OO_X^{\oplus \binom{s}{1}}(-c+1)
\longrightarrow\cdots\\
& \longrightarrow 
\OO_X^{\oplus \binom{s}{c-1}}(-1)
\longrightarrow 
\OO_X^{\oplus \binom{s}{c}}
\longrightarrow 
{\textstyle\bigwedge^c \mathcal{Q}_X}
\longrightarrow 0.
\end{aligned}$$
Since 
$\big({\textstyle\bigwedge^b \mathcal{Q}_X}\big)^\vee
= \big({\textstyle\bigwedge^{s-1-b} \mathcal{Q}_X}\big)(-1)$, 
we have an exact sequence
$$
\begin{aligned}
0&\longrightarrow \OO_X(b-s)
\longrightarrow \OO_X^{\oplus \binom{s}{s-1}}(b-s+1)
\longrightarrow\cdots\\
& \longrightarrow 
\OO_X^{\oplus \binom{s}{b+2}}(-2)
\longrightarrow 
\OO_X^{\oplus \binom{s}{b+1}}(-1)
\longrightarrow 
\big({\textstyle\bigwedge^b \mathcal{Q}_X}\big)^{\vee}
\longrightarrow 0.
\end{aligned}$$
So we can find a quasi-isomorphism 
$\mathcal{C}^\bullet\to 
\EEnd (\mathcal{T}_X)
=\mathcal{T}_X\otimes \mathcal{T}_X^\vee
$ with $\mathcal{C}^{-p} = \OO_X(-p-1)^{\oplus(\cdots)}$ for $p\geq 0$. 
There is a spectral sequence 
$$E_1^{-p,q}(-j)
= 
H^q(X,\mathcal{C}^{-p}(-j))
\Longrightarrow 
H^{q-p}(X,\EEnd_X(\mathcal{T}_X)(-j)).
$$
It suffices to show 
$$E_1^{-p,q}(-j) = H^q(X,\OO_X(-p-1-j))^{\oplus (\cdots)} = 0$$
for $q-p=k>j\geq 0$ and $p\geq 0$. 
In this case, it could be nonzero only when $q=s-1$, but it is zero since $-p-j-1\geq -q=1-s>-s$. 

By the above exact sequence and induction, ${\bigwedge^{0} \mathcal{Q}_X}
\oplus 
{\bigwedge^{1} \mathcal{Q}_X}
\oplus \cdots \oplus 
{\bigwedge^{d} \mathcal{Q}_X}
$ generates 
$\OO_X\oplus \OO_X(-1)\oplus \ldots\oplus \OO_X(-d)$ for any $0\leq d\leq s-1$. 
By \cite[\href{https://stacks.math.columbia.edu/tag/0A9V}{0A9V}]{Stacks}, for $d=s-1$, they generate $D_{qc}(X)$. 
\end{proof}

\begin{proof}[Proof of Theorem \ref{thm:Tistilting}]
Since $\mathcal{T}_Y$ is a vector bundle, it suffices to show that $\mathcal{T}_Y$ is a generator and 
$H^{>0}(Y,\EEnd(\mathcal{T}_Y))=0$. 

Note that 
$$\EEnd_Y(\mathcal{T}_Y) = 
\EEnd_X(\mathcal{T}_X)|_{Y}.$$
Tensoring the Koszul complex \eqref{eq:KoszulfornHil} with 
$\EEnd_X(\mathcal{T}_X)$, we get a quasi-isomorphism 
$\mathcal{C}^\bullet\to 
\EEnd_Y(\mathcal{T}_Y)$ with 
$$\mathcal{C}^{-p}
= 
\EEnd_X(\mathcal{T}_X)\otimes \OO_X(-p)^{\oplus (\cdots)}
= 
\EEnd_X(\mathcal{T}_X)(-p)^{\oplus (\cdots)}
.$$
There is a spectral sequence 
$$E_1^{-p,q} = H^q(X,\mathcal{C}^{-p})
\Longrightarrow H^{q-p}(Y,\EEnd(\mathcal{T}_Y)).$$
It suffices to show 
$$E_1^{-p,q}=H^q(X,\EEnd_X(\mathcal{T}_X)(-p))^{\oplus (\cdots)}=0$$
for $p<q$. 
This is exactly Lemma \ref{lem:Lambdaistitling}. 

Adjunction gives
$$
R\Hom_Y(\mathcal{T}_Y,\mathcal{K})
=
R\Hom_X(\mathcal{T}_X,\mathcal{K}),$$
where we identify $\mathcal{K}$ with its pushforward along the closed inclusion $Y\to X$. Since $\mathcal{T}_X$ is a generator for $X$, $\mathcal{T}_Y$ is a generator for $Y$. 
\end{proof}

The following lemma is a consequence of tilting theory. 

\begin{lemma}\label{lem:inadd}
Let $\mathcal T$ be a tilting bundle over a quasi-compact, quasi-separated scheme $\mathfrak X$. 
If a vector bundle $\mathcal{F}$ satisfies 
$$
\Ext_{\mathfrak{X}}^{>0}(\mathcal{T},\mathcal{F})
=0=
\Ext_{\mathfrak{X}}^{>0}(\mathcal{F},\mathcal{T}),$$
then $\mathcal{F}\in \operatorname{add}(\mathcal{T})$.
\end{lemma}
\begin{proof}
Let $A = \End_{\mathfrak{X}}(\mathcal{T})$. 
Derived Morita theory gives a triangulated equivalence 
$$
R\Hom_{\mathfrak{X}}(\mathcal{T},-): D_{qc}(\mathfrak{X})
\stackrel{\sim}\longrightarrow 
D(A): \mathcal{T}\otimes_A -. $$
Then the first vanishing implies 
$R\Hom_{\mathfrak{X}}(\mathcal{T},\mathcal{F})$ is an $A$-module, say $M$. 
The second vanishing implies 
$$\Ext_A^{>0}(M,A)=0.$$
Since $\mathcal{F}$ is compact, the module $M$ is also compact. 
By \cite[\href{https://stacks.math.columbia.edu/tag/07LT}{07LT}]{Stacks}, 
$M$ is perfect, i.e. it is quasi-isomorphic to a bounded complex $P^\bullet$ of finitely generated projective $A$-modules. 
By splitting the surjective differential of $P^{\bullet}$ at the highest degree, we can assume $P^{>0}=0$.
Now $\Hom_A(P^\bullet,A)$ is the complex computing $\Ext^\bullet_A(M,A)$. 
By the second vanishing, the perfect complex $\Hom_A(P^\bullet,A)$ is quasi-isomorphic to 
$\Hom_A(M,A)$. 
By splitting the surjective differential and using the fact that 
$$P^\bullet = \Hom_{A^{op}}(\Hom_A(P^\bullet,A),A^{op}),$$
we can conclude that $M$ itself is projective, hence belongs to $\operatorname{add}(A)$. 
As a result, $\mathcal{F} = \mathcal{T}\otimes_A M\in \operatorname{add}(\mathcal{T})$.  
\end{proof}

\subsection{Results of Theorem \ref{thm:Section3}}

The diagram \eqref{eq:mixed-main-diagram} gives the following corollary of Theorem \ref{thm:Section3}. 

\begin{corollary}
\label{cor:mixed}
We have 
\begin{equation}\label{eq:intro-mixed}
\begin{gathered}
R^{>0}p_{n,*}q_{n-1}^*(\mathcal{P}_{n-1}\otimes\mathcal{B}_{n-1}^{\otimes b})=0,\\
\mathcal{P}_n\otimes\mathcal{B}_n^{\otimes b}
\twoheadrightarrow p_{n,*}q_{n-1}^*(\mathcal{P}_{n-1}\otimes\mathcal{B}_{n-1}^{\otimes b}).
\end{gathered}
\end{equation}
Since $Y_{n,0}^-=Y_{n,0}$, we have 
$\mathcal{P}_n\simeq p_{n*}q_{n-1}^*\mathcal{P}_{n-1}$. 
\end{corollary}
\begin{proof}
Let $\rho^-_{n,b}:Y^-_{n,b}\to H_n$ be the restriction of $\rho_{n,b}$.
Finite base change in \eqref{eq:mixed-main-diagram} and
Theorem~\ref{thm:Section3} give
\[
Rp_{n,*}q_{n-1}^*(\mathcal P_{n-1}\otimes\mathcal B_{n-1}^{\otimes b})
\simeq \rho^-_{n,b,*}\mathcal O_{Y^-_{n,b}}.
\]
The right-hand side is concentrated in degree zero.
The closed immersion $Y^-_{n,b}\hookrightarrow Y_{n,b}$ and exactness
of finite pushforward give the required surjection.
For $b=0$, the two schemes coincide, so this surjection is an isomorphism.
\end{proof}

Recall the line bundle $\mathcal{L}$ defined in \eqref{eq:tau-N-Hn}.

\begin{proposition}
\label{prop:mixed}
We have 
$$R^{>0}\alpha_*\big(q_{n-1}^*(\mathcal{P}_{n-1}\otimes{\textstyle\bigwedge^b\mathcal{B}_{n-1}})\otimes \mathcal{L}\big)=0.$$
\end{proposition}
\begin{proof}
Since $p_n=\pi_n\circ \alpha$ and $\pi_{n,*}$ is exact and faithful, it suffices to show 
\begin{equation}\label{eq:equivmixed}
R^{>0}p_{n,*}\big(q_{n-1}^*(\mathcal{P}_{n-1}\otimes{\textstyle\bigwedge^b\mathcal{B}_{n-1}})\otimes \mathcal{L}\big)=0.
\end{equation}
The sequence \eqref{eq:tau-N-Hn} gives 
\begin{equation}\label{eq:sesinmixed}
0\textstyle\to 
q_{n-1}^*\bigwedge^b\mathcal{B}_{n-1}\otimes \mathcal{L}
\to 
p_n^*\bigwedge^{b+1}\mathcal{B}_{n}
\to 
q_{n-1}^*
\bigwedge^{b+1}\mathcal{B}_{n-1}\to 0
.
\end{equation}
Note that by the projection formula and Corollary \ref{cor:mixed}, we have 
\begin{align*}
\textstyle 
Rp_{n,*}\big(
q_{n-1}^*\mathcal{P}_{n-1}\otimes 
p_n^*\bigwedge^{b+1}\mathcal{B}_{n}
\big)
&\textstyle 
= 
\mathcal{P}_{n}\otimes 
\bigwedge^{b+1}\mathcal{B}_{n}. 
\end{align*}
By Corollary \ref{cor:mixed}, we have 
$$\begin{gathered}
\textstyle 
R^{>0}p_{n,*}q_{n-1}^*\big(
\mathcal{P}_{n-1}\otimes 
\bigwedge^{b+1}\mathcal{B}_{n-1}
\big)=0,\\
\textstyle 
\mathcal{P}_{n}\otimes 
\bigwedge^{b+1}\mathcal{B}_{n}
\twoheadrightarrow
p_{n,*}q_{n-1}^*\big(
\mathcal{P}_{n-1}\otimes 
\bigwedge^{b+1}\mathcal{B}_{n-1}
\big).
\end{gathered}$$
Tensoring \eqref{eq:sesinmixed} with $q_{n-1}^*\mathcal{P}_{n-1}$ and applying $p_{n,*}$ gives the vanishing \eqref{eq:equivmixed} by the long exact sequence. 
\end{proof}

\begin{corollary}\label{coro:Ralpha=0}
We have 
$$
\begin{gathered}
\textstyle 
R^{>0}\alpha_{*}(
q_{n-1}^*\mathcal{P}_{n-1}\otimes 
q_{n-1}^*
(\bigwedge^b \mathcal{B}_{n-1})^\vee)=0,\\
\textstyle 
R^{>0}\alpha_{*}(
q_{n-1}^*\mathcal{P}_{n-1}^\vee\otimes 
q_{n-1}^*
\bigwedge^b \mathcal{B}_{n-1})=0.
\end{gathered}$$
\end{corollary}
\begin{proof}
As in the proof of Proposition \ref{prop:mixed}, it suffices to prove the $R^{>0}p_{n,*}$ version. 
The sequence \eqref{eq:tau-N-Hn} gives 
$$p_{n}^*\det\mathcal{B}_n = 
q_{n-1}^*\det\mathcal{B}_{n-1}\otimes \mathcal{L}.$$
So we have 
$$
\begin{aligned}
\textstyle
q_{n-1}^*
(\bigwedge^b \mathcal{B}_{n-1})^\vee
& = 
\textstyle
q_{n-1}^*
(\bigwedge^{n-1-b} \mathcal{B}_{n-1}\otimes \det \mathcal{B}_{n-1}^\vee)\\
& = 
\textstyle
q_{n-1}^*
(\bigwedge^{n-1-b} \mathcal{B}_{n-1})\otimes \mathcal{L}\otimes 
p_n^* \det\mathcal{B}_n^{\vee}.
\end{aligned}
$$
By the projection formula, the first vanishing follows from Proposition \ref{prop:mixed}. 

By Haiman's duality formula \cite[Proposition 3.7.4]{Haiman}, 
\begin{equation}
\label{eq:Haimanduality}
\mathcal{P}_m^\vee \simeq 
\mathcal{P}_m\otimes (\det\mathcal{B}_m)^{\vee} 
\otimes \operatorname{sgn}_{S_m}.
\end{equation}
Thus, non-equivariantly, we have 
$$\textstyle 
q_{n-1}^*\mathcal{P}_{n-1}^\vee
\simeq 
q_{n-1}^*\mathcal{P}_{n-1}
\otimes 
\mathcal{L}\otimes 
p_n^* \det\mathcal{B}_n^{\vee}.$$
By the projection formula, the second vanishing also follows from Proposition \ref{prop:mixed}. 
\end{proof}

\subsection{Local structure of $\alpha$}
Now let us return to the local structure of $\alpha$. 

\begin{lemma}
\label{lem:localsplit}
Over $Y$, we have 
$$q_{n-1}^* \mathcal{B}_{n-1}|_Y 
\simeq \mathcal{Q}|_Y\oplus \mathcal{O}_Y^{n-s}. $$
\end{lemma}
\begin{proof}
Note that we can embed 
$$H_{n-1,n}\stackrel{\subset}\longrightarrow
\mathbb{P}_{U_n}(\pi_n^*\mathcal{B}_n). 
$$
This embedding can be made compatible with \eqref{eq:nHil-is-proj}. 
Recall 
$$H_{n-1,n}=\mathbb{P}_{U_n}^\vee(\kappa_{U_n})=\mathbb{P}_{U_n}^\vee(\kappa_{\pi_n})$$
for the relative dualizing sheaf $\kappa_{\pi_n}$; see Remark \ref{rmk:relativeshv}. 
From its construction, we have a natural counit
$$\pi_n^*\mathcal{B}_n^\vee \twoheadrightarrow \kappa_{\pi_n}. $$
This induces the closed embedding 
$\mathbb{P}_{U_n}^\vee(\kappa_{\pi_n})\subset\mathbb{P}_{U_n}^\vee(\pi_n^*\mathcal{B}_n^\vee)$. 
Let $S = \mathcal{O}_{H_n\times \mathbb{A}^2,\xi}$ and $S/\mathfrak{a} = \mathcal{O}_{U_n,\xi}$ as in Section \ref{sec:locstrofalpha}. 
Recall the surjective map $S^{\oplus s}\to \kappa_{S/\mathfrak{a}}$. 
By trivializing $\pi_n^*\mathcal{B}_n^\vee$, 
we can lift the counit to $C$ in the following diagram:
$$\xymatrix@C=4pc{
& S^{\oplus n}\ar@{->>}[d]
\ar@{..>}[dl]_{C}\\
S^{\oplus s}
\ar[r]&\kappa_{S/\mathfrak{a}}\ar[r]&0}$$
Since the horizontal morphism induces an isomorphism on fibers at $\xi$, a maximal minor of $C$ is a unit in $S$. Thus, after a basis change, we can assume $C=(I_s,0)$. So we have 
$$\xymatrix{
0\ar[r] & 
\mathcal{L}\ar[r]\ar@{=}[d]&
p_{n}^*\mathcal{B}_{n}\ar[r]
\ar@{<-}[d]_{\left[\substack{I_s\\0}\right]}&
q_{n-1}^*\mathcal{B}_{n-1}\ar[r]\ar@{<..}[d]&0\\
0\ar[r] & 
\mathcal{O}_Y(-1)\ar[r] & S^{\oplus s}\ar[r] & \mathcal{Q}|_Y\ar[r] & 0.
}$$
Taking the cokernel gives $q_{n-1}^* \mathcal{B}_{n-1}|_Y 
\simeq \mathcal{Q}|_Y\oplus \mathcal{O}_Y^{n-s}$. 
\end{proof}

\subsection{Proof of Theorem \ref{thm:Section4}}
By Corollary \ref{coro:Ralpha=0}, 
$$
\begin{aligned}
\textstyle
\Ext_Y^{>0}(q_{n-1}^*\mathcal{P}_{n-1}|_Y\,,\,\,
q_{n-1}^*\bigwedge^b \mathcal{B}_{n-1}|_Y)&= 0,\\
\textstyle
\Ext_Y^{>0}(
q_{n-1}^*\bigwedge^b \mathcal{B}_{n-1}|_Y\,,\,\,
q_{n-1}^*\mathcal{P}_{n-1}|_Y)
& =0.
\end{aligned}
$$
By Lemma \ref{lem:localsplit}, we get 
$$
\begin{aligned}
\textstyle
\Ext_Y^{>0}(q_{n-1}^*\mathcal{P}_{n-1}|_Y\,,\,\,
\bigwedge^b \mathcal{Q}|_Y)&= 0,\\
\textstyle
\Ext_Y^{>0}(
\bigwedge^b \mathcal{Q}|_Y\,,\,\,
q_{n-1}^*\mathcal{P}_{n-1}|_Y)
& =0.
\end{aligned}$$
As a result, $q_{n-1}^*\mathcal{P}_{n-1}|_Y$ satisfies the conditions of Lemma \ref{lem:inadd}.
In particular, 
$$q_{n-1}^*\mathcal{P}_{n-1}|_Y\in \operatorname{add}(\mathcal{T}_Y). $$

\section{A cohomological model for the intersections}\label{sec:regular}

Let $n\geq2$ and $k\geq1$ be integers, and let 
$\lambda\vdash n$. Let
$\mu(1),\ldots,\mu(k)\vdash n-1$ be distinct partitions with $|\lambda/\mu(i)|=1$. 
Consider 
\begin{equation}\label{eq:V}
V = \operatorname{span}(x^ay^b: 
(a,b)=\lambda/\mu(i)
\text{ for some $1\leq i\leq k$}).
\end{equation}
Note that we can identify $V$ with a subspace of $\operatorname{soc}(\mathbb{C}[x,y]/I_\lambda)_0$.
Consider the following diagram
$$
\xymatrix{
\mathbb{P}(V)\ar[d]_{f}\ar[r]^-{\tilde{f}}
& H_{n-1,n}\ar[dl]^{q_{n-1}}\ar[dr]_{p_{n}}\ar[r]^-{\alpha} & U_n\ar[d]^{\pi_n}\\
H_{n-1}&&H_n,}$$
where the morphism $f$ is given by 
$f(L)=I_\lambda+L$ while 
$\tilde{f}(L)=(I_\lambda,f(L))$. 
Both morphisms are closed embeddings. 
Moreover, the image of $\tilde{f}$ is contained in $\alpha^{-1}(\xi)$ for $\xi=(I_\lambda,0)$.

\begin{maintheorem}\label{thm:Section5}
(1) 
The higher cohomology groups of $f^*\mathcal{P}_{n-1}^\vee$ vanish:
$$H^{>0}(\mathbb{P}(V),f^*\mathcal{P}_{n-1}^\vee)=0.$$
(2) 
For the zeroth cohomology, we have
\begin{equation}\label{eq:H0isinter}
H^0(\mathbb{P}(V),f^*\mathcal{P}_{n-1}^\vee)^\vee
\simeq \frac{\mathbb{C}[x_1,y_1,\ldots,x_{n-1},y_{n-1}]}{J_{\mu(1)}+\cdots+J_{\mu(k)}}.
\end{equation}
\end{maintheorem}

Let us explain the construction of \eqref{eq:H0isinter}. 
Since $X_n\subset H_n\times (\mathbb{A}^2)^n$ is closed, 
we have a surjective morphism
$$\mathcal{O}_{H_n}[x_1,y_1,\ldots,x_n,y_n]\twoheadrightarrow \mathcal{P}_{n},$$
whose fiber at $I\in H_n$ is the coordinate quotient
$\mathbb{C}[x_1,y_1,\ldots,x_n,y_n]\twoheadrightarrow\mathcal{P}_n(I)$. 
Similarly,
$$\mathcal{O}_{\mathbb{P}(V)}[x_1,y_1,\ldots,x_{n-1},y_{n-1}]\twoheadrightarrow f^*\mathcal{P}_{n-1}. $$
It induces 
$$\mathbb{C}[x_1,y_1,\ldots,x_{n-1},y_{n-1}]
\longrightarrow
H^0(\mathbb{P}(V),f^*\mathcal{P}_{n-1})
\longrightarrow 
H^0(\mathbb{P}(V),f^*\mathcal{P}_{n-1}^\vee)^\vee. 
$$
We will show that this map is surjective with kernel $J_{\mu(1)}+\cdots+J_{\mu(k)}$.

\subsection{Castelnuovo--Mumford regularity}

In this subsection, we prove the main tool for proving Theorem \ref{thm:Section5}.

\begin{remark}\label{rem:CMregu}
Recall Castelnuovo--Mumford regularity
\cite[\href{https://stacks.math.columbia.edu/tag/089X}{089X}]{Stacks}. 
A coherent sheaf $\mathcal{E}$ on projective space $\mathbb{P}^n$ is called $m$-regular if 
$$H^{i}(\mathbb{P}^n,\mathcal{E}(m-i))=0\text{ for all }i>0. $$
For example, a line bundle $\OO(d)$ is $m$-regular if and only if $d\geq -m$. 
Clearly, $\mathcal{E}$ is $m$-regular if and only if $\mathcal{E}(m)$ is 0-regular. 

We recall some properties of $0$-regular sheaves. 
\begin{itemize}
    \item 
A 0-regular sheaf is globally generated; see  \cite[\href{https://stacks.math.columbia.edu/tag/08A8}{08A8}]{Stacks}. 

\item 
The restriction of a 0-regular vector bundle $\mathcal{E}$ to any hyperplane $H$ (and hence any linear subspace) is still 0-regular; 
see \cite[\href{https://stacks.math.columbia.edu/tag/08A5}{08A5}]{Stacks}. 
In fact, this can be seen from the short exact sequence
$0\to \mathcal{E}(-1)\to \mathcal{E}\to \mathcal{E}|_H\to 0$. 

\item If $\mathcal{E}$ is 0-regular, then so are $\mathcal{E}(1),\mathcal{E}(2),\ldots$; see 
\cite[\href{https://stacks.math.columbia.edu/tag/08A6}{08A6}]{Stacks}. 
In particular, 
$$H^{>0}(\mathbb{P}^n,\mathcal{E}(-1))=0.$$
\end{itemize}
\end{remark}

\begin{theorem}\label{thm:0regularandsum}
Let $\mathcal{E}$ be a 0-regular vector bundle on $\mathbb{P}^m$. 
For any surjective map $H\to H^0(\mathbb{P}^m,\mathcal{E})$, 
let
$$K_i = \ker(H\to 
\mathcal{E}([e_i])),$$
where $e_0,e_1,\ldots,e_{m}$ form the standard basis of $\mathbb{C}^{m+1}$. 
Then we have an isomorphism 
$$H/(K_0+\cdots+K_m)
\simeq H^0(\mathbb{P}^m,\mathcal{E}^\vee)^\vee,$$
induced by 
$H\to 
H^0(\mathbb{P}^m,\mathcal{E})
\to
H^0(\mathbb{P}^m,\mathcal{E}^\vee)^\vee$.
\end{theorem}

\begin{proof}
The kernel of $H\to H^0(\mathbb P^m,\mathcal E)$ is contained
in every $K_i$. Passing to the quotient, we may assume
$H=H^0(\mathbb P^m,\mathcal E)$.
Define the vector bundle $\mathcal{K}$ on $\mathbb{P}^m$ by  
$$
0\longrightarrow \mathcal{K}(-1)\longrightarrow 
H\otimes \OO_{\mathbb{P}^m}
\longrightarrow 
\mathcal{E}\longrightarrow 0.
$$
Then the first map in the above short exact sequence is equivalent to a linear map
$\gamma:K\otimes \mathbb{C}^{m+1}\rightarrow H$ 
for $K=H^0(\mathbb{P}^m,\mathcal{K})$. 
By definition, 
$\gamma(K\otimes\mathbb C e_i)=K_i$ and the image of $\gamma$ is $K_0+\cdots+K_m$. 

The long exact cohomology sequence shows that $\mathcal{K}$ is $0$-regular, hence globally generated.
We thus get a locally split exact sequence
$$
K\otimes \OO_{\mathbb{P}^m}(-1)\longrightarrow 
H\otimes \OO_{\mathbb{P}^m}
\longrightarrow 
\mathcal{E}\longrightarrow 0.
$$
Dualizing the sequence, we get 
\begin{equation}\label{eq:beforetakinggs}
0\longrightarrow \mathcal{E}^\vee \longrightarrow 
H^\vee\otimes \OO_{\mathbb{P}^m}
\longrightarrow 
K^\vee\otimes \OO_{\mathbb{P}^m}(1). 
\end{equation}
Taking global sections, we arrive at
$$0\longrightarrow 
H^0(\mathbb{P}^m,\mathcal{E}^\vee)
\longrightarrow H^\vee \longrightarrow K^\vee\otimes (\mathbb{C}^{m+1})^\vee. 
$$
Dualizing again, we get 
$$
K\otimes \mathbb{C}^{m+1}
\stackrel{\gamma}\longrightarrow
H
\longrightarrow 
H^0(\mathbb{P}^m,\mathcal{E}^\vee)^\vee\longrightarrow 0,$$
and this concludes the proof.
\end{proof}

\begin{corollary}
\label{coro:0regularandsum}
We can identify $\mathbb{P}^d\subset \mathbb{P}^m$ when $0\leq d<m$. 
Under the hypotheses of Theorem \ref{thm:0regularandsum}, we have 
$$H/(K_0+\cdots+K_d)
\simeq H^0(\mathbb{P}^d,\mathcal{E}^\vee|_{\mathbb{P}^d})^\vee.$$
\end{corollary}
\begin{proof}
The proof is the same as the proof of Theorem \ref{thm:0regularandsum}, except that we take global sections of \eqref{eq:beforetakinggs} over $\mathbb{P}^d$ instead. 
\end{proof}

\begin{lemma}
\label{lem:quotis0reg}
Over $\mathbb{P}^m$, the sheaf $\mathcal{T}:=
{\textstyle \bigwedge^\bullet 
\mathcal{Q}}$
is $0$-regular, where $\mathcal{Q}$ is the universal quotient bundle. 
\end{lemma}
\begin{proof}
We prove the equivalent vanishing statement 
$$H^{i}(\mathbb{P}^{m},\mathcal{T}(-j))=0\qquad \text{for $i>0$ and $0\leq j\leq i$}. $$
The proof is the same as the proof of Lemma \ref{lem:Lambdaistitling}. 
In this case, we omit the resolution of $\mathcal{T}^\vee$, so there is one fewer twist. 
\end{proof}

\subsection{Proof of Theorem \ref{thm:Section5} (1)}
Let us adopt the notation \eqref{eq:XYZ}. 
We have a tautological sequence on $\mathbb{P}(V)$:
\begin{equation}\label{eq:tauseqoverP}
0\longrightarrow 
\mathcal{O}_{\mathbb{P}(V)}(-1)
\longrightarrow 
(\mathbb C[x,y]/I_\lambda)\otimes \mathcal{O}_{\mathbb{P}(V)}
\longrightarrow 
f^*\mathcal{B}_{n-1}
\longrightarrow 
0.
\end{equation}
It gives 
$\mathcal{O}_{\mathbb{P}(V)}(-1)\simeq f^* \det \mathcal{B}_{n-1}^\vee$. 
By Haiman's duality formula \cite[Proposition 3.7.4]{Haiman}, 
(non-equivariantly)
\begin{equation}
\label{eq:HaimandualityoverP}
f^*\mathcal{P}_{n-1}(-1)
\simeq f^*\mathcal{P}_{n-1}\otimes f^* \det \mathcal{B}_{n-1}^\vee\simeq 
(f^*\mathcal{P}_{n-1})^{\vee}. 
\end{equation}
As a result, it suffices to prove $f^*\mathcal{P}_{n-1}$ is $0$-regular; see Remark \ref{rem:CMregu}. 

Let $\xi=(I_\lambda,0)\in U_n$ and $s=\dim_{\mathbb C}\operatorname{soc}(\mathbb{C}[x,y]/I_\lambda)_0$. 
Note that the closed embedding
$\tilde{f}$ 
identifies $\mathbb{P}(V)$ with a
linear subspace of $\alpha^{-1}(\xi)=Z\cong \mathbb{P}^{s-1}$.
It suffices to show 
$q_{n-1}^*\mathcal{P}_{n-1}|_{Z}$
is $0$-regular; see Remark \ref{rem:CMregu}. 
By Theorem \ref{thm:Section4}, 
$$q_{n-1}^*\mathcal{P}_{n-1}|_{Z}\in \operatorname{add}(\mathcal{T}_Z).$$
By Lemma \ref{lem:quotis0reg}, it is $0$-regular.
Theorem \ref{thm:Section5} (1) follows.

\subsection{Proof of Theorem \ref{thm:Section5} (2)}
Since $X_n$ is closed in $H_n\times (\mathbb{A}^2)^n$, we have a natural surjective morphism
$$\mathcal{O}_{H_n}[x_1,y_1,\ldots,x_n,y_n]\twoheadrightarrow \mathcal{P}_n.$$
Put $\mathcal{M}=\alpha_*q_{n-1}^*\mathcal{P}_{n-1}$. 
By Corollary \ref{cor:mixed}, 
$\pi_{n,*}\mathcal{M}=\mathcal{P}_n$. The counit for the finite map $\pi_n$ is surjective, so we get 
$\pi_n^*\mathcal{P}_n\twoheadrightarrow 
\mathcal{M}$. 
Hence,
$$\mathcal{O}_{U_n}[x_1,y_1,\ldots,x_{n},y_{n}]
\twoheadrightarrow 
\pi_n^*\mathcal{P}_n
\twoheadrightarrow 
\mathcal{M}. $$
Since $x_n,y_n$ act through $\mathcal{O}_{U_n}$, we obtain $$\mathcal{O}_{U_n}[x_1,y_1,\ldots,x_{n-1},y_{n-1}]
\twoheadrightarrow 
\mathcal{M}. $$
Taking the fiber at $\xi$ and using flat base change, we obtain a surjective map 
$$\mathbb{C}[x_1,y_1,\ldots,x_{n-1},y_{n-1}]
\twoheadrightarrow
H^0(Y,q_{n-1}^*\mathcal{P}_{n-1}|_Y)\otimes S/\mathfrak{m}
.$$
Here we used the notation \eqref{eq:Sam}. 

Moreover, the restriction
$$
H^0(Y,q_{n-1}^*\mathcal{P}_{n-1}|_Y)
\to 
H^0(Z,q_{n-1}^*\mathcal{P}_{n-1}|_Z)$$ is surjective. This can be seen as follows.
Let us adopt the notation from \eqref{eq:tilting}. 
By Theorem \ref{thm:Section4}, this will follow from the surjectivity of 
$$
H^0(Y,\mathcal{T}_Y)
\to 
H^0(Z,\mathcal{T}_Z),
$$
which holds since the composition 
$$
H^0(X,\mathcal{T}_X)
\to 
H^0(Y,\mathcal{T}_Y)
\to 
H^0(Z,\mathcal{T}_Z)$$
is surjective. 

Hence, we get a surjective map
\begin{equation}\label{eq:surjweneed}
\mathbb{C}[x_1,y_1,\ldots,x_{n-1},y_{n-1}]
\twoheadrightarrow
H^0(Z,q_{n-1}^*\mathcal{P}_{n-1}|_Z).
\end{equation}

In the proof of Theorem \ref{thm:Section5} (1), we showed that $q_{n-1}^*\mathcal{P}_{n-1}|_Z$ is $0$-regular on $Z\cong \mathbb{P}^{s-1}$. Hence, Theorem \ref{thm:Section5} (2) follows from
Theorem \ref{thm:0regularandsum} and Corollary \ref{coro:0regularandsum} by setting $$H:=\mathbb{C}[x_1,y_1,\ldots,x_{n-1},y_{n-1}],\text{ and } \mathcal{E}=q_{n-1}^*\mathcal{P}_{n-1}|_Z.$$  Notice that by the above construction, $K_{i-1}=J_{\mu(i)}$.

\section{The intersection formula via localization and primitive modules}
\label{sec:localization}

\subsection{Intersection formulas}
Recall that $T=(\mathbb{C}^\times)^2$ acts on $H_n$ naturally, and $S_n$ acts trivially on $H_n$. 
Note that $V$ in \eqref{eq:V} is a $T$-representation, and $f:\mathbb{P}(V)\to H_{n-1}$ is $T$-equivariant. 
Since $x^ay^b$ has $T$-weight $q^at^b$, the weights of $V$ are distinct, so the $T$-fixed points are 
$$\mathbb{P}(V)^T=\{[v_i]:
i=1,\ldots,k\},\qquad 
v_i=x^{a_i}y^{b_i},
(a_i,b_i)=\lambda/\mu(i).$$
Recall that  for a partition $\lambda$, 
$$T_\lambda:=\prod_{(i,j)\in \lambda} q^it^j = q^{n(\lambda')}t^{n(\lambda)}.$$
It is the $T$-weight of $\det\mathcal{B}_n$ at $I_\lambda\in H_n$. Therefore, the tangent space $T_{[v_i]}\mathbb{P}(V)=\Hom(\mathbb{C}v_i,V/\mathbb{C}v_i)$ has $T$-weights 
\[\{q^{a_j-a_i}t^{b_j-b_i}=T_{\mu(i)}/T_{\mu(j)}\mid 1\leq j\neq i\leq k\}.\]

\begin{maintheorem}\label{thm:Section6}
We have 
$$
\begin{aligned}
&\operatorname{Frob}_{q,t}\left(H^0(\mathbb{P}(V),f^*\mathcal{P}_{n-1}^{\vee})^\vee\right)
=
\operatorname{Frob}_{q,t}
\big(\mathsf{GH}_{\mu(1)}\cap \cdots \cap \mathsf{GH}_{\mu(k)}\big)
\\
&\qquad = \sum_{i=1}^k 
\left(\prod_{j\neq i}
\frac{T_{\mu(j)}}{T_{\mu(j)}-T_{\mu(i)}}\right)
\widetilde{H}_{\mu(i)}[X;q,t].
\end{aligned}
$$
\end{maintheorem}
\begin{proof}
The first equality follows from Theorem \ref{thm:Section5}. 
Indeed, 
$$
\begin{aligned}
H^0(\mathbb{P}(V),f^*\mathcal{P}_{n-1}^{\vee})
& \simeq 
(J_{\mu(1)}+\cdots+J_{\mu(k)})^\perp \\
& = 
J_{\mu(1)}^\perp\cap \cdots\cap J_{\mu(k)}^\perp \\
& = \mathsf{GH}_{\mu(1)}\cap \cdots \cap \mathsf{GH}_{\mu(k)}. 
\end{aligned}$$
This isomorphism reverses the grading, by the definition of the pairing in Section \ref{2-11}. 
Since $S_{n-1}$-representations are self-dual, the first equality follows. 

The second equality follows from the Atiyah--Bott localization formula for the equivariant K-theory of $\mathbb{P}(V)$ \cite{CG97}.
Writing $\mathcal E=f^*\mathcal P_{n-1}$, Theorem \ref{thm:Section5} (1) and the fixed-point formula give
\begin{align*}
\operatorname{Frob}_{q,t}H^0(\mathbb P(V),\mathcal E^\vee)
&=\sum_{i=1}^k\frac{\widetilde H_{\mu(i)}[X;q^{-1},t^{-1}]}{\prod_{j\ne i}(1-T_{\mu(j)}/T_{\mu(i)})}\\
&=\sum_{i=1}^k\left(\prod_{j\ne i}\frac{T_{\mu(j)}^{-1}}{T_{\mu(j)}^{-1}-T_{\mu(i)}^{-1}}\right)\widetilde H_{\mu(i)}[X;q^{-1},t^{-1}].
\end{align*}
Taking ordinary duals, hence substituting $q^{-1},t^{-1}$ for $q,t$, gives
\begin{align*}
    \operatorname{Frob}_{q,t}\left(H^0(\mathbb{P}(V),f^*\mathcal{P}_{n-1}^{\vee})^\vee\right)=&\operatorname{Frob}_{q^{-1},t^{-1}}H^0(\mathbb P(V),\mathcal E^\vee)\\
    =&\sum_{i=1}^k\left(\prod_{j\ne i}\frac{T_{\mu(j)}}{T_{\mu(j)}-T_{\mu(i)}}\right)\widetilde H_{\mu(i)}[X;q,t].
\end{align*}
This is the asserted formula.
\end{proof}

\subsection{Primitive Garsia--Haiman modules}
\label{sec:compl}
Let us give an extension of Lemma \ref{lem:quotis0reg}. 
We  first recall the definition of an distributive lattice. 
A \emph{lattice} is a partially ordered set $(L,\leq)$ in which every pair $x,y\in L$ has a greatest lower bound $x\wedge y$ and a least upper bound $x\vee y$. 
We also write $(L,\wedge,\vee)$ for a lattice. 
The lattice is \emph{distributive} if, for  $x,y,z\in L$,
$$
\begin{aligned}
x\wedge(y\vee z)&=(x\wedge y)\vee(x\wedge z),\\
x\vee(y\wedge z)&=(x\vee y)\wedge(x\vee z).
\end{aligned}$$
For example, the subsets of a given set form a distributive lattice under intersection and union; 
the subspaces of a vector space form a (generally not distributive) lattice under intersection and sum.

\begin{lemma}
\label{lem:Kiasso}
Assume $\mathcal{E}\in \operatorname{add}(\mathcal{T})$ for $\mathcal{T}$ as in Lemma \ref{lem:quotis0reg}. 
Let $K_i\subset H$ be as in Theorem \ref{thm:0regularandsum}. 
Then $K_0,\ldots,K_{m}$ generate a distributive lattice.
\end{lemma}
\begin{proof}
As in the proof of Theorem \ref{thm:0regularandsum}, we can assume $H=H^0(\mathbb{P}^m,\mathcal{T})$. 
Note that 
the assertion for $\mathcal{E}_1$ and $\mathcal{E}_2$ is equivalent to the assertion for $\mathcal{E}_1\oplus \mathcal{E}_2$.
Thus it suffices to prove the assertion when $\mathcal{E}=\mathcal{T}$. 

The same complex in the proof of Lemma \ref{lem:Lambdaistitling} gives 
$H^0(\mathbb{P}^m,{\textstyle\bigwedge^\bullet}\mathcal{Q})
={\textstyle\bigwedge^\bullet}\mathbb{C}^{m+1}$.
Thus, each
$$K_i = \ker\left[
{\textstyle\bigwedge^\bullet}\mathbb{C}^{m+1}\to {\textstyle\bigwedge^\bullet}(\mathbb{C}^{m+1}/\mathbb{C}e_i)
\right] = e_i\wedge {\textstyle\bigwedge^{\bullet}}\mathbb{C}^{m+1}$$
is a coordinate space of ${\textstyle\bigwedge^{\bullet}}\mathbb{C}^{m+1}$. Therefore, $K_0,\ldots,K_m$ generate a distributive lattice. 
\end{proof}

Using Lemma \ref{lem:Kiasso} and Theorem \ref{thm:Section4}, we can derive the following distributive lattice property of Garsia--Haiman modules. 

\begin{theorem}\label{thm:GHmodisass}
The Garsia--Haiman modules
$\mathsf{GH}_{\mu(1)},\ldots,\mathsf{GH}_{\mu(k)}$ generate a distributive lattice.
\end{theorem}
\begin{proof}
Put $R=\mathbb{C}[x_1,y_1,\ldots,x_{n-1},y_{n-1}]$, and let
$Z=\alpha^{-1}(I_\lambda,0)$ and
$\mathcal{E}=q_{n-1}^*\mathcal{P}_{n-1}|_Z$. Choose coordinates on $Z$ by extending the monomials
$v_1,\ldots,v_k$ to a basis of
$\operatorname{soc}(\mathbb{C}[x,y]/I_\lambda)_0$.
By Theorem \ref{thm:Section4},
$\mathcal{E}\in\operatorname{add}(\mathcal{T}_Z)$, and
\eqref{eq:surjweneed} gives a natural surjection
$R\twoheadrightarrow H^0(Z,\mathcal{E})$.
For each $1\leq i\leq k$, the composite evaluation map
$$
R\longrightarrow H^0(Z,\mathcal{E})
\longrightarrow \mathcal{E}([v_i])
=\mathcal{P}_{n-1}(I_{\mu(i)})
$$
is the coordinate quotient with kernel $J_{\mu(i)}$.
Thus Lemma \ref{lem:Kiasso}, applied with $H=R$, shows that
$J_{\mu(1)},\ldots,J_{\mu(k)}$ generate a distributive lattice.

The pairing in Section \ref{2-11} is perfect in each bidegree.
Consequently, for bihomogeneous ideals $I,J\subset R$ of finite codimension,
$$
(I+J)^\perp=I^\perp\cap J^\perp,
\qquad
(I\cap J)^\perp=I^\perp+J^\perp,
\qquad
(I^\perp)^\perp=I.
$$
Every ideal in the lattice generated by the $J_{\mu(i)}$ is
bihomogeneous and of finite codimension. Hence taking orthogonal
complements gives an order-reversing lattice isomorphism from this
lattice to the lattice generated by
$J_{\mu(i)}^\perp=\mathsf{GH}_{\mu(i)}$.
Since the dual of a distributive lattice is distributive, the result follows.
\end{proof}

For any nonempty subset $A\subseteq [k]$, recall that 
$\mathsf{GH}_A = 
\bigcap_{i\in A}\mathsf{GH}_{\mu(i)}$.

\begin{maintheorem}
\label{lem:primitivepro}
There exists a submodule $\mathsf{GH}'_A$ of $\mathsf{GH}_A$ for each nonempty subset $A\subseteq [k]$, called a 
\emph{primitive Garsia--Haiman module}, such that  
$$
\mathsf{GH}_{\mu(1)}
+\cdots+\mathsf{GH}_{\mu(k)}
=\bigoplus_{A\neq \varnothing} \mathsf{GH}'_A,
\qquad 
\mathsf{GH}_{A}=\bigoplus_{B\supseteq A} \mathsf{GH}'_B.$$    
\end{maintheorem}
\begin{proof}
We pick any $S_{n-1}$-equivariant complement $\mathsf{GH}'_A$ 
of $\sum_{B\supsetneqq A}\mathsf{GH}_B$ inside $\mathsf{GH}_A$. 
The fact that 
$$\mathsf{GH}_{\mu(1)}
+\cdots+\mathsf{GH}_{\mu(k)}
=\sum_{A\neq \varnothing} \mathsf{GH}'_A,\qquad 
\mathsf{GH}_{A}=\sum_{B\supseteq A} \mathsf{GH}'_B$$
follows directly from the definition of $\mathsf{GH}'_A$. Hence, it suffices to show the sum is direct. 
After ordering the summands by increasing cardinality, this is equivalent to the vanishing of
\begin{equation}\label{eq:disum}
\mathsf{GH}'_{B}\cap \sum_{C\neq B,|C|\geq |B|} \mathsf{GH}'_C.
\end{equation}
Notice that
\begin{align*}
\eqref{eq:disum} 
& \subseteq 
\mathsf{GH}_{B}\cap 
\sum_{\substack{|C|\geq |B|\\C\neq B}}
\mathsf{GH}_C
\stackrel{(*)}=\sum_{\substack{|C|\geq |B|\\C\neq B}}
\mathsf{GH}_{B}\cap \mathsf{GH}_C
=
\sum_{\substack{|C|\geq |B|\\C\neq B}}
\mathsf{GH}_{B\cup C}
\subseteq 
\sum_{D\supsetneqq B}
\mathsf{GH}_{D}.
\end{align*}
Here $(*)$ uses Theorem \ref{thm:GHmodisass}. 
Since 
$\mathsf{GH}'_{B}$ is a complement to
$\sum_{D\supsetneqq B}
\mathsf{GH}_{D}$, 
\eqref{eq:disum} is zero. 
\end{proof}

\begin{corollary}\label{cor:conseq-dimensions}
For $\varnothing\ne A\subseteq[k]$ of cardinality $a$, one has
\begin{equation*}\label{eq:conseq-primitive-dim}
\dim_{\mathbb C}\mathsf{GH}'_A
=\frac{N!(a-1)!(k-a)!}{k!}
=\frac{N!}{k {k-1 \choose a-1}}.
\end{equation*}
Moreover,
\begin{equation*}\label{eq:conseq-sum-dim}
\dim_{\mathbb C}\left(\sum_{i=1}^k\mathsf{GH}_{\mu(i)}\right)
=N!\sum_{a=1}^k\frac1a.
\end{equation*}
\end{corollary}

\begin{proof}
By the M\"obius inversion, Theorem~\ref{lem:primitivepro}, and Corollary~\ref{cor:N!/k}, we get
\[
\begin{aligned}
\dim_{\mathbb C}\mathsf{GH}'_A
&=N!\sum_{j=0}^{k-a}(-1)^j\binom{k-a}{j}\frac1{a+j}\\
&=N!\int_0^1u^{a-1}(1-u)^{k-a}\,du
=\frac{N!(a-1)!(k-a)!}{k!}.
\end{aligned}
\]
The second identity follows from the first one and 
Theorem~\ref{lem:primitivepro}.
\end{proof}

Next, let us give a cohomology model for the primitive Garsia--Haiman modules. 

\begin{theorem}
\label{thm:Section6prim}
Introduce $T_A = \prod_{i\in A}T_{\mu(i)}$ for any $A\subseteq [k]$.
For any nonempty subset $A\subseteq [k]$ of cardinality $a>0$, we have
$$
\begin{aligned}
& T_A\cdot 
\operatorname{Frob}_{q,t}\left(H^{a-1}\big(\mathbb{P}(V),
f^*(\mathcal{P}_{n-1}\otimes
\det\mathcal{B}_{n-1}^{\vee,{\otimes a}})\big)
\right)
=
\operatorname{Frob}_{q,t}
\mathsf{GH}'_A\\
&
\qquad = -T_A\cdot 
\sum_{i=1}^k 
\frac{1}{\prod_{j\neq i}(1-T_{\mu(j)}/T_{\mu(i)})}
\cdot 
\frac{\widetilde{H}_{\mu(i)}[X;q,t]}{(-T_{\mu(i)})^{a}}.
\end{aligned}
$$
\end{theorem}
\begin{proof}
By Theorem \ref{thm:Section6}, 
\begin{align*}
\operatorname{Frob}_{q,t}\mathsf{GH}_A
& = \sum_{i\in A}\prod_{\substack{\alpha\in A\\\alpha\neq i}}
\frac{T_{\mu(\alpha)}}{T_{\mu(\alpha)}-T_{\mu(i)}}
\widetilde{H}_{\mu(i)}[X;q,t]\\
& = - T_A
\sum_{i\in A}
\prod_{\substack{\alpha\in A\\\alpha\neq i}}
\frac{1}{1-T_{\mu(\alpha)}/T_{\mu(i)}}\cdot 
\frac{\widetilde{H}_{\mu(i)}[X;q,t]}{(-T_{\mu(i)})^{a}}\\
& = -T_A
\sum_{i=1}^k
\frac{\prod_{j\in [k]\setminus A}(1-T_{\mu(j)}/T_{\mu(i)})}
{\prod_{j\neq i}(1-T_{\mu(j)}/T_{\mu(i)})}\cdot 
\frac{\widetilde{H}_{\mu(i)}[X;q,t]}{(-T_{\mu(i)})^a}\\
& = -T_A
\sum_{i=1}^k
\frac{\sum_{B\supseteq A}
\prod_{\beta\in B\setminus A}(-T_{\mu(\beta)}/T_{\mu(i)})}
{\prod_{j\neq i}(1-T_{\mu(j)}/T_{\mu(i)})}\cdot 
\frac{\widetilde{H}_{\mu(i)}[X;q,t]}{(-T_{\mu(i)})^a}\\
& = 
-\sum_{B\supseteq A}T_B
\sum_{i=1}^k
\frac{1}
{\prod_{j\neq i}(1-T_{\mu(j)}/T_{\mu(i)})}\cdot 
\frac{\widetilde{H}_{\mu(i)}[X;q,t]}{(-T_{\mu(i)})^{|B|}}.
\end{align*}
By Theorem \ref{lem:primitivepro}, we have 
$$\operatorname{Frob}_{q,t}\mathsf{GH}'_{B}
= -T_B\cdot 
\sum_{i=1}^k
\frac{1}
{\prod_{j\neq i}(1-T_{\mu(j)}/T_{\mu(i)})}\cdot 
\frac{\widetilde{H}_{\mu(i)}[X;q,t]}{(-T_{\mu(i)})^{|B|}}.$$
This proves the second equality. 

To prove the first equality, we need the following property. 
For $0<a\leq k$, we have 
$$
H^{p}(\mathbb{P}(V),f^*\mathcal{P}_{n-1}(-a))
=0\text{ unless }p=a-1.$$
When $a\leq p$, this follows from the fact that $f^*\mathcal{P}_{n-1}$ is $0$-regular, and hence so is 
$f^*\mathcal{P}_{n-1}(p-a)$; see the proof of Theorem \ref{thm:Section5} (1) and Remark \ref{rem:CMregu}.
Now assume $0\leq p\leq a-2$. 
By Serre duality
$$
H^{p}(\mathbb{P}(V),f^*\mathcal{P}_{n-1}(-a))
\simeq 
H^{k-1-p}(\mathbb{P}(V),f^*\mathcal{P}^\vee_{n-1}(a-k))^\vee.
$$
By Haiman's duality formula \cite[Proposition 3.7.4]{Haiman} (see \eqref{eq:HaimandualityoverP}), it is isomorphic to 
$$H^{k-1-p}(\mathbb{P}(V),f^*\mathcal{P}_{n-1}(a-k-1))^\vee.$$
Since $(a-k-1)+(k-1-p)=a-p-2\geq 0$ and $k-1-p\geq k-1-(a-2) = k-a+1>0$, 
the vanishing follows from the fact that $f^*\mathcal{P}_{n-1}(a-p-2)$ is $0$-regular, see Remark \ref{rem:CMregu}.

Taking determinants in \eqref{eq:tauseqoverP}, 
we get a non-equivariant isomorphism $$f^*(\det\mathcal{B}_{n-1}^\vee)
\simeq\mathcal{O}_{\mathbb{P}(V)}(-1).$$ 
Then the first formula follows from these facts and the Atiyah--Bott formula.
\end{proof}

Theorem \ref{thm:Section6prim} shows that
\[
T_A^{-1}\cdot\operatorname{Frob}_{q,t}\mathsf{GH}'_A
\]
depends only on $|A|$. Primitive modules indexed by subsets of the same size therefore have the same bigraded representation structure after a grading shift, although this does not itself supply canonical identifications between the chosen complements.

The following theorem should be viewed as the analogue of Theorem \ref{thm:Section6prim} for $A=\varnothing$. 

\begin{theorem}
\label{thm:Section6prim2}
We have 
$$\begin{aligned}
& \operatorname{Frob}_{q,t}\big(H^0(\mathbb{P}(V),f^*\mathcal{P}_{n-1})\big)=
\operatorname{Frob}_{q,t}
\big(\mathsf{GH}_{\mu(1)}+\cdots +\mathsf{GH}_{\mu(k)}\big)
\\
&\qquad = \sum_{i=1}^k 
\left(\prod_{j\neq i}
\frac{T_{\mu(i)}}{T_{\mu(i)}-T_{\mu(j)}}\right)
\widetilde{H}_{\mu(i)}[X;q,t].
\end{aligned}$$
Note that this is distinct from the formula in Theorem \ref{thm:Section6}. 
\end{theorem}
\begin{proof}
By Theorem \ref{lem:primitivepro} and Theorem \ref{thm:Section6prim}, 
\begin{align*}
\operatorname{Frob}_{q,t}
\left(
\sum_{i=1}^k \mathsf{GH}_{\mu(i)}
\right)
& = 
-\sum_{A\neq \varnothing}
T_A\cdot 
\sum_{i=1}^k 
\frac{1}{\prod_{j\neq i}(1-T_{\mu(j)}/T_{\mu(i)})}
\cdot 
\frac{\widetilde{H}_{\mu(i)}[X;q,t]}{(-T_{\mu(i)})^{a}}\\
& \hspace{-5pc} 
= -\sum_{i=1}^k 
\frac{1}{\prod_{j\neq i}(1-T_{\mu(j)}/T_{\mu(i)})}
\left(\sum_{A\neq \varnothing}
(-1)^a\frac{T_A}{T_{\mu(i)}^a}
\right)\cdot \widetilde{H}_{\mu(i)}[X;q,t]\\
& \hspace{-5pc} 
= -\sum_{i=1}^k 
\frac{1}{\prod_{j\neq i}(1-T_{\mu(j)}/T_{\mu(i)})}
\left(\prod_{j=1}^k(1-T_{\mu(j)}/T_{\mu(i)})-1
\right)\cdot \widetilde{H}_{\mu(i)}[X;q,t]\\
& \hspace{-5pc} 
= \sum_{i=1}^k 
\frac{1}{\prod_{j\neq i}(1-T_{\mu(j)}/T_{\mu(i)})}
\cdot \widetilde{H}_{\mu(i)}[X;q,t].
\end{align*}
This is the second equality. 
Since $f^*\mathcal{P}_{n-1}$ is $0$-regular, $H^{>0}(\mathbb{P}(V),f^*\mathcal{P}_{n-1})=0$; see Remark \ref{rem:CMregu}. 
By Atiyah--Bott localization, the first equality holds. 
\end{proof}

\begin{remark}
Recall that the nabla operator on the symmetric functions with coefficients in $\mathbb{Q}(q,t)$ is the $q,t$-linear map satisfying $\nabla(\widetilde{H}_\lambda)= T_\lambda\cdot \widetilde{H}_\lambda$. It was originally introduced in  Bergeron--Garsia's Science Fiction \cite{BG}, and played a  ubiquitous role in  the field of Macdonald polynomials \cite{BGHT}; see, for example,   $\nabla e_n$  in the diagonal coinvariants \cite{Hai-02, HHLRA, CM} and  $\nabla  s_\lambda$  in the Loehr--Warrington conjecture  \cite{LW0,BMMPS}.
\begin{itemize}
    
\item 
 
Theorem \ref{thm:Section6prim} and Theorem \ref{thm:Section6prim2} implies
$$\operatorname{Frob}_{q,t}\mathsf{GH}'_A = (-1)^{a+1}T_A\cdot \nabla^{-a} \left(\operatorname{Frob}_{q,t}
\big(\mathsf{GH}_{\mu(1)}+\cdots+\mathsf{GH}_{\mu(k)}\big)\right).
$$

\item 
Applying Serre duality to Theorem \ref{thm:Section6prim}, or computing directly, one can get 
$$\begin{aligned}
\operatorname{Frob}_{q,t}\mathsf{GH}_{A}'=& T_{[k]\setminus A}^{-1}\cdot\operatorname{Frob}_{q,t}
\left(H^{k-a}(\mathbb{P}(V),f^*(\mathcal{P}_{n-1}^{\vee}\otimes \det\mathcal{B}_{n-1}^{\vee,\otimes(k-a)}))^\vee\right)\\
=& T_{[k]\setminus A}^{-1}\cdot
\sum_{i=1}^k 
\left(\prod_{j\neq i}\frac{T_{\mu(j)}}{T_{\mu(j)}-T_{\mu(i)}}\right)
(-T_{\mu(i)})^{k-a}\cdot \widetilde{H}_{\mu(i)}[X;q,t]\\
=&T_{[k]\setminus A}^{-1}\cdot(-\nabla)^{k-a}\left(\operatorname{Frob}_{q,t}
\big(\mathsf{GH}_{\mu(1)}\cap \cdots \cap \mathsf{GH}_{\mu(k)}\big)\right).
\end{aligned}$$
Here the second formula follows from the Atiyah--Bott localization formula while the last one follows from Theorem \ref{thm:Section6}. In particular, if $A=[k]$, we 
get Theorem \ref{thm:Section6}.

\item 
Combining the above two, we get
\begin{equation*}\label{eq:conseq-sum-nabla}
\operatorname{Frob}_{q,t}
\big(\mathsf{GH}_{\mu(1)}+ \cdots + \mathsf{GH}_{\mu(k)}\big)
=\frac{(-1)^{k-1}}{T_{[k]}}
\nabla^k\left(\operatorname{Frob}_{q,t}
\big(\mathsf{GH}_{\mu(1)}\cap \cdots \cap \mathsf{GH}_{\mu(k)}\big)\right).
\end{equation*}
This can also be obtained directly from Theorem~\ref{thm:Section6} and Theorem~\ref{thm:Section6prim2}.

    \item 
By Theorem \ref{thm:Section6prim2}, it is not hard to see 
\begin{equation}\label{eq:intersect}
    H^0(\mathbb{P}(V),f^*\mathcal{P}_{n-1})
\simeq \frac{\mathbb{C}[x_1,y_1,\ldots,x_{n-1},y_{n-1}]}{J_{\mu(1)}\cap \cdots\cap J_{\mu(k)}}.
\end{equation}
Compare this with Theorem \ref{thm:Section5}. 
In fact, by the surjectivity of \eqref{eq:surjweneed} and the surjectivity of the restriction (see Remark \ref{rem:CMregu})
\[H^0(Z,q_{n-1}^*\mathcal{P}_{n-1}|_Z)\rightarrow H^0(\mathbb{P}(V),f^*\mathcal{P}_{n-1}),\] 
we get
$$H^0(\mathbb{P}(V),f^*\mathcal{P}_{n-1})
\simeq \frac{\mathbb{C}[x_1,y_1,\ldots,x_{n-1},y_{n-1}]}{\bigcap_{P\in \mathbb{P}(V)}J_{f(P)}}\twoheadrightarrow \frac{\mathbb{C}[x_1,y_1,\ldots,x_{n-1},y_{n-1}]}{J_{\mu(1)}\cap \cdots\cap J_{\mu(k)}},$$
where $J_{I}$ is the kernel of $\mathbb{C}[x_1,\ldots,y_{n-1}]\to \mathcal{P}_{n-1}(I)$. 

The pairing in Section \ref{2-11} gives
$(\bigcap_{i=1}^kJ_{\mu(i)})^\perp
=\sum_{i=1}^k\mathsf{GH}_{\mu(i)}$.  Hence, \eqref{eq:intersect} holds as they have the same dimension by Theorem \ref{thm:Section6prim2}.
\end{itemize}
\end{remark}

\end{document}